# The Lingering of Gradients: Theory and Applications

Zeyuan Allen-Zhu
zeyuan@csail.mit.edu
Microsoft Research AI

David Simchi-Levi
dslevi@mit.edu
MIT

Xinshang Wang
xinshang@mit.edu
MIT

(version 2)*

**Abstract**

Classically, the time complexity of a first-order method is estimated by its number of gradient computations. In this paper, we study a more refined complexity by taking into account the "lingering" of gradients: once a gradient is computed at $x_k$, the additional time to compute gradients at $x_{k+1}, x_{k+2}, \ldots$ may be reduced.

We show how this improves the running time of several first-order methods. For instance, if the "additional time" scales linearly with respect to the traveled distance, then the "convergence rate" of gradient descent can be improved from $1/T$ to $\exp(-T^{1/3})$. On the application side, we solve a hypothetical revenue management problem on the Yahoo! Front Page Today Module with 4.6m users to $10^{-6}$ error using only 6 passes of the dataset; and solve a real-life support vector machine problem to an accuracy that is two orders of magnitude better comparing to the state-of-the-art algorithm.

## 1 Introduction

First-order methods play a fundamental role in large-scale machine learning and optimization tasks. In most scenarios, the performance of a first-order method is represented by its *convergence rate*: the relationship between the optimization error $\varepsilon$ and the number of gradient computations $T$. This is meaningful because in most applications, the time complexities for evaluating gradients at different points are of the same magnitude. In other words, the worse-case time complexities of first-order methods are usually proportional to a fixed parameter times $T$.

In certain large-scale settings, if we have already spent time computing the (full) gradient at $x$, perhaps we can use such information to reduce the time complexity to compute full gradients at other points near $x$. We call this the "lingering" of gradients, because the gradient at $x$ may be partially reused for future consideration, but will eventually fade away once we are far from $x$.

In this paper, we consider an important class of optimization problems in which algorithms can exploit the lingering of gradients and thus converge faster. Formally, consider the (finite-sum) stochastic convex minimization problem:

$$\min_{x \in \mathbb{R}^d} \left\{ f(x) \stackrel{\text{def}}{=} \tfrac{1}{n} \sum_{i=1}^n f_i(x) \right\} \ . \tag{1.1}$$

Then, could it be possible that *whenever* $x$ is sufficiently close to $y$, for at least a large fraction of indices $i \in [n]$, we have $\nabla f_i(x) \approx \nabla f_i(y)$? In other words, if $\nabla f_1(x), \ldots, \nabla f_n(x)$ are already calculated at some point $x$, can we reuse a large fraction of them to approximate $\nabla f(y)$?

---

*An extended abstract of a preliminary version of this paper has appeared in the conference NeurIPS 2018 with title "The Lingering of Gradients: How to Reuse Gradients Over Time", and the current paper is a significant extension to it.

**Example 1.** In the problem of matching customers to resources, $f_i(x)$ represents the marginal profit of the $i$-th customer under bid-price vector $x \in \mathbb{R}_+^d$ over $d$ items. In many applications (see Section 2.1), $\nabla f_i(x)$ only depends on customer $i$'s preferences under $x$.

If the bid-price vector $x \in \mathbb{R}_+^d$ changes by a small amount to $y$, then for a large fraction of customers $i$, their most profitable items may not change, and thus $\nabla f_i(x) \approx \nabla f_i(y)$. Indeed, imagine if one of the items is Xbox, and its price drops by 5%, perhaps 90% of the customers will not change their minds about buying or not. We shall demonstrate this using real-life data.

**Example 2.** In classification problems, $f_i(x)$ represents the loss value for "how well training sample $i$ is classified under predictor $x$". For any sample $i$ that has a large margin under predictor $x$, its gradient $\nabla f_i(x)$ may stay close to $\nabla f_i(y)$ whenever $x$ is close to $y$.

Formally, let $f_i(x) = \max\{0, 1 - \langle x, a_i \rangle\}$ be the hinge loss (or its smoothed variant if needed) with respect to the $i$-th sample $a_i \in \mathbb{R}^d$. If the margin $|1 - \langle x, a_i \rangle|$ is sufficiently large, then moving from $x$ to a nearby point $y$ should not affect the sign of $1 - \langle x, a_i \rangle$, and thus not change the gradient. Therefore, if samples $a_1, \ldots, a_n$ are sufficiently diverse, then a large fraction of them should incur large margins and have the same gradients when $x$ changes by little.

## 1.1 Summary of Main Results and Contributions

We assume in this paper that, given any point $x \in \mathbb{R}^d$ and index $i \in [n]$, one can efficiently evaluate a "lingering radius" $\delta(x, i)$. The radius satisfies the condition that for every point $y$ that is within distance $\delta(x, i)$ from $x$, the stochastic gradient $\nabla f_i(y)$ is equal to $\nabla f_i(x)$. We make two remarks:

- We use "equal to" for the purpose of proving theoretical results. In practice and in our experiments, it suffices to use approximate equality such as $\|\nabla f_i(x) - \nabla f_i(y)\| \leq 10^{-10}$.

- By "efficient" we mean $\delta(x, i)$ is computable in the same complexity as evaluating $\nabla f_i(x)$. This is reasonable because when $\nabla f_i(x)$ is an explicit function of $x$, it is usually easy to tell how sensitive it is to the input $x$. (We shall include such examples in our experiments.)

If we denote by $B(x, r)$ the set of indices $j$ satisfying $\delta(x, j) < r$, and if we travel to some point $y$ that is at most distance $r$ from $x$, then we only need to re-evaluate the (stochastic) gradients $\nabla f_j(y)$ for $j \in B(x, r)$. Intuitively, one should expect $|B(x, r)|$ to grow as a function of $r$ if the data points are sufficiently diverse.

**Better Convergence Rate in Theory.** To present the simplest theoretical result, we modify gradient descent (GD) to take into account the lingering of gradients. At a high level, we run GD, but during its execution, we maintain a decomposition of the indices $\Lambda_0 \cup \cdots \cup \Lambda_t = \{1, 2, \ldots, n\}$ where $t$ is logarithmic in $n$. Now, whenever we need $\nabla f_i(x_k)$ for some $i \in \Lambda_p$, we approximate it by $\nabla f_i(x_{k'})$ for a point $k'$ that was visited at most $2^p$ steps ago. Our algorithm makes sure that such $\nabla f_i(x_{k'})$ is available in memory.

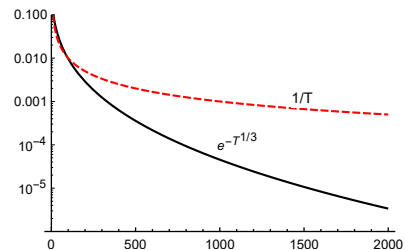

Figure 1: $e^{-T^{1/3}}$ vs $1/T$

We prove that the performance of our algorithm depends on how $|B(x, r)|$ grows in $r$. Formally, let $T$ be the total number of stochastic gradient computations divided by $n$, and suppose $|B(x, r)| \leq O(r^\beta)$. Then, our algorithm finds a point $x$ with $f(x) - f(x^*) \leq \widetilde{O}(T^{-1/(1-\beta)})$ if $\beta \in (0, 1)$, or $f(x) - f(x^*) \leq 2^{-\Omega(T^{1/3})}$ if $\beta = 1$. In contrast, traditional GD satisfies $f(x) - f(x^*) \leq O(T^{-1})$.

**Faster Algorithm in Practice.** We also design an algorithm that practically maximizes the use of gradient lingering. We take the SVRG method [21, 38] as the prototype because it is widely applied in large-scale settings. Recall that SVRG uses gradient estimator $\nabla f(\widetilde{x}) - \nabla f_i(\widetilde{x}) + \nabla f_i(x_k)$ to estimate the full gradient $\nabla f(x_k)$, where $\widetilde{x}$ is the so-called snapshot point (which was visited at most $n$ steps ago) and $i$ is a random index. At a high level, we modify SVRG so that the index $i$ is only generated from those whose stochastic gradients need to be recomputed, and ignore those such that $\nabla f_i(x_k) = \nabla f_i(\widetilde{x})$. This can further reduce the variance of the gradient estimator, and improve the running time.

**Application to packing LPs.** Our algorithms serve as tools for solving a variety of packing linear programs (LPs), including those widely used by revenue-maximization policies [13, 34]. In this paper, we solve a packing LP of this form on the Yahoo! Front Page Today Module application [7, 24] with 4.6 million users to $10^{-6}$ error (or $10^{-12}$ dual error) using only 6 passes of the dataset.

**Application to SVM.** Our algorithms also apply to training support vector machine (SVM), one of the most classical supervised learning model for classification tasks. On the Adult dataset of LibSVM [12], we manage to minimize the SVM training objective to $10^{-5}$ error in 30 passes of the dataset. In contrast, PEGASOS, arguably the most popular method for SVM [33], cannot minimize this objective even to $10^{-3}$ error within 90 passes.

## 1.2 Related Work

**Variance Reduction.** The SVRG method was independently proposed by Johnson and Zhang [21], Zhang et al. [38], and belong to the class of stochastic methods using the so-called variance-reduction technique [5, 9, 21, 25, 29–32, 37, 38]. The common idea behind these methods is to use some full gradient of the past to approximate future, but they do not distinguish which $\nabla f_i(x)$ can "linger longer in time" among all indices $i \in [n]$ for different $x$.

Arguably the two most widely applied variance-reduction methods are SVRG and SAGA [9]. They have complementary performance depending on the internal structural of the dataset [6], so we compare to both in our experiments.

A practical modification of SVRG is to use an approximate full gradient (as opposed to the exact full gradient) of the past to approximate future. This is studied by [17, 22, 23], and we refer to this method as SCSG due to [22, 23].

**Reuse Gradients.** Some researchers have exploited the internal structure of the dataset to speed up first-order methods. That is, they use $\nabla f_i(x)$ to approximate $\nabla f_j(x)$ when the two data samples $i$ and $j$ are sufficiently close. This is orthogonal to our setting because we use $\nabla f_i(x)$ to approximate $\nabla f_i(y)$ when $x$ and $y$ are sufficiently close. In the extreme case when all the data samples are identical, they have $\nabla f_i(x) = \nabla f_j(x)$ for every $i, j$ and thus stochastic gradient methods converge as fast as full gradient ones. For this problem, Hofmann et al. [18] introduce a variant of SAGA, Allen-Zhu et al. [6] introduce a variant of SVRG and a variant of accelerated coordinate descent.

Other authors study how to reduce gradient computations at the snapshot points of SVRG [17, 22]. This is also orthogonal to the idea of this paper, and can be added to our algorithms for even better performance (see Section 4.2).

**A Preliminary Version.** An extended abstract of a preliminary version of this paper has appeared in the conference NeurIPS 2018, and the current paper is a significant extension to it. Specifically, the current version has three more major contributions.

- First, we now provide theories for a more general assumption on the lingering radius (the current Assumption 2 allows $\beta \in (0, 1]$ while the conference version only allows $\beta = 1$).

- Second, we now apply our methods also to the task of support vector machines (Section 6).

- Third, we now provide theories showing that the assumption of lingering radius indeed holds when data is sufficiently random (Section 7).

Besides these major contributions, we have additionally applied our technique to the SCSG method and conducted more thorough experiments.

## 1.3 Roadmap

In Section 2, we introduce notations for this paper and give setups for our packing LP and SVM applications. In Section 3, we prove our main theoretical result on the improved convergence rate for gradient descent under the aforementioned assumption $|B(x,r)| \leq O(r^\beta)$. In Section 4, we introduce our practical algorithm by incorporating the lingering of gradients into SVRG and SCSG. Using real-life datasets, we apply our algorithms to packing LP in Section 5 and to SVM in Section 6. Finally, in Section 7, we provide theoretical support for the assumption $|B(x,r)| \leq O(r^\beta)$ using randomness of the data.

# 2 Notions and Problem Formulation

We denote by $\|\cdot\|$ the Euclidean norm, and $\|\cdot\|_\infty$ the infinity norm. Recall the notion of Lipschitz smoothness (it has other equivalent definitions, see textbook [27]).

**Definition 2.1.** *A function $f \colon \mathbb{R}^d \to \mathbb{R}$ is L-Lipschitz smooth (or L-smooth for short) if*
$$\forall x, y \in \mathbb{R}^d \colon \|\nabla f(x) - \nabla f(y)\| \leq L\|x - y\| \enspace.$$

We propose the following model to capture the lingering of gradients.

**Definition 2.2.** *For every $x \in \mathbb{R}^d$ and index $i \in [n]$, let $\delta(x,i) \geq 0$ be the* lingering radius *of $\nabla f_i(x)$, meaning that*[1]
$$\nabla f_i(x) = \nabla f_i(y) \text{ for all } y \in \mathbb{R}^d \text{ with } \|y - x\| \leq \delta(x,i)$$

*Accordingly, for every $r \geq 0$ we use $B(x,r)$ to denote the set of indices $j$ satisfying $\delta(x,j) < r$:*
$$B(x,r) \stackrel{\text{def}}{=} \{j \in [n] \,|\, \delta(x,j) < r\} \enspace.$$

In other words, as long as we travel within distance $\delta(x,i)$ from $x$, the gradient $\nabla f_i(x)$ can be reused to represent $\nabla f_i(y)$. Our main assumption of this paper is that

**Assumption 1.** *Each $\delta(x,i)$ can be computed in the same time complexity as $\nabla f_i(x)$.*

Under Assumption 1, if at some point $x$ we have already computed $\nabla f_i(x)$ for all $i \in [n]$, then we can compute $\delta(x,i)$ as well within the same time complexity for every $i \in [n]$, and sort the indices $i \in [n]$ in increasing order of $\delta(x,i)$. In the future, if we arrive at any point $y$, we can calculate $r = \|x - y\|$ and use
$$\nabla' = \tfrac{1}{n} \Big( \sum_{i \notin B(x,r)} \nabla f_i(x) + \sum_{i \in B(x,r)} \nabla f_i(y) \Big)$$
to represent $\nabla f(y)$. The time to compute $\nabla'$ is only proportional to $|B(x,r)|$.

**Definition 2.3.** *We denote by $T_{\text{time}}$ the gradient complexity, which equals how many times $\nabla f_i(x)$ and $\delta(x,i)$ are calculated, divided by $n$.*

---
[1]Recall that, in practice, one should replace the exact equality with, for instance, $\|\nabla f_i(x) - \nabla f_i(y)\| \leq 10^{-10}$. To present the simplest statements, we do not introduce such an extra parameter.

In computing $\nabla'$ above, the gradient complexity is $|B(x,r)|/n$. If we always set $\delta(x,i) = 0$ then $|B(x,r)| = n$ and the gradient complexity for computing $\nabla'$ remains 1. However, if the underlying Problem (1.1) is nice enough so that $|B(x,r)|$ becomes an increasing function of $r$, then the gradient complexity for computing $\nabla'$ can be less than 1. We can thus hope for designing faster algorithms.

## 2.1 Packing Linear Program

Consider the LP relaxation of a canonical revenue management problem in which a manager needs to sell $d$ different resources to $n$ customers. Let $b_j \geq 0$ be the capacity of resource $j \in [d]$; let $p_{i,j} \in [0,1]$ be the probability that customer $i \in [n]$ will purchase a unit of resource $j$ if offered resource $j$; and let $r_j$ be the revenue for each unit of resource $j$. We want to offer each customer one and only one candidate resource, and let $y_{i,j}$ be the probability we offer customer $i$ resource $j$. The following is the standard LP relaxation for this problem:[2]

$$\max_{y \geq 0} \sum_{i \in [n], j \in [d]} r_j p_{i,j} y_{i,j}$$
$$\text{s.t.} \sum_{i \in [n]} p_{i,j} y_{i,j} \leq b_j, \quad \forall j \in [d] \quad (2.1)$$
$$\sum_{j \in [d]} y_{i,j} = 1, \quad \forall i \in [n].$$

This LP (2.1) and its variants have repeatedly found many applications, including adwords/ad allocation problems [3, 11, 13, 16, 19, 26, 36, 39], and revenue management for airline and service industries [8, 15, 20, 28, 34, 35]. Some authors also study the online version of solving LPs [1, 2, 10, 14].

A standard way to reduce (2.1) to convex optimization is by regularization, see for instance Zhong et al. [39]. Let us subtract the maximization objective by a regularizer

$$R(y) \stackrel{\text{def}}{=} \mu \sum_{i \in [n]} \overline{p_i} \sum_{j \in [d]} y_{i,j} \log y_{i,j},$$

where $\overline{p_i} \stackrel{\text{def}}{=} \max_{i \in [n]} p_{i,j}$ and $\mu > 0$ is some small regularization weight. Then, after transforming to the dual, we have a new minimization problem

$$\min_{x \geq 0} \mu \sum_{i=1}^{n} \overline{p_i} \cdot \log Z_i + \sum_{j=1}^{d} x_j b_j, \quad (2.2)$$

where

$$Z_i = \sum_{j=1}^{d} \exp\left(\frac{(r_j - x_j) p_{i,j}}{\overline{p_i} \mu}\right).$$

If we let $f_i(x) \stackrel{\text{def}}{=} \mu n \overline{p_i} \cdot \log Z_i + \langle x, b \rangle$, then (2.2) reduces to Problem (1.1). We conduct empirical studies on this packing LP problem in Section 5.

*Remark* 2.4. Any solution $x$ (usually known as the *bid price*) to (2.2) naturally gives back a solution $y$ for the primal (2.1), by setting

$$y_{i,j} = \exp\left(\frac{(r_j - x_j) p_{i,j}}{\overline{p_i} \mu}\right) / Z_i. \quad (2.3)$$

---
[2]The constraint $\sum_{j \in [d]} y_{i,j} = 1$ here can be replaced with any other positive constant without loss of generality.

## 2.2 Support Vector Machine

Classifying data is one of the most foundational tasks in machine learning. Suppose we are given data points $a_1, \ldots, a_n \in \mathbb{R}^d$ each belonging to one of two classes. We use $b_i = 1$ to denote that data point $i$ belongs to the first class, and $b_i = -1$ to denote that data point $i$ belongs to the second.

The (soft-margin) support vector machine task is to minimize the following objective

$$f(x) = \frac{\lambda}{2}\|x\|^2 + \frac{1}{n}\sum_{i=1}^{n} \max\{0, 1 - b_i \langle x, a_i \rangle\} , \qquad (2.4)$$

where $\lambda$ is the weight of the regularizer which encourages the objective to find a solution with large classification margin. If we set $f_i(x) = \frac{\lambda}{2}\|x\|^2 + \max\{0, 1 - b_i\langle x, a_i\rangle\}$, then (2.4) reduces to Problem (1.1).

In this formulation, the SVM objective $f(x)$ is not Lipschitz smooth, making some of the popular practical methods unable to apply (at least in theory). For such reason, people also study the smoothed version of SVM as follows.[3]

$$f^\mu(x) = \frac{1}{n}\sum_{i=1}^{n} f_i^\mu(x) \quad \text{where} \quad f_i^\mu(x) = \frac{\lambda}{2}\|x\|^2 + \begin{cases} 0, & b_i\langle x, a_i\rangle \geq 1; \\ 1 - \frac{\mu}{2} - b_i\langle x, a_i\rangle, & b_i\langle x, a_i\rangle \leq 1 - \mu; \\ \frac{1}{2\mu}(1 - b_i\langle x, a_i\rangle)^2, & \text{otherwise.} \end{cases} \qquad (2.5)$$

Above, $\mu \geq 0$ is a smoothing parameter. The larger $\mu$ is, the more Lipschitz smooth the objective $f^\mu(x)$ becomes. We conduct empirical studies on this SVM problem in Section 6.

## 3 GD with Lingering Radius

In this section, we consider a convex function $f(x) = \frac{1}{n}\sum_{i=1}^{n} f_i(x)$ that is $L$-smooth. Recall from textbooks (e.g., [27]) that if gradient descent (GD) is applied for $T$ iterations, starting at $x_0 \in \mathbb{R}^d$, then we can arrive at a point $x$ with $f(x) - f(x^*) \leq O\big(\frac{\|x_0 - x^*\|^2}{T}\big)$. This is the $\frac{1}{T}$ convergence rate.

To improve on this theoretical rate, we make the following assumption on $B(x, r)$:

**Assumption 2.** *There exists $\alpha \in [0, 1], \beta \in (0, 1], C > 0$ such that,*

$$\forall x \in \mathbb{R}^d, r \geq 0: \quad \frac{|B(x, r)|}{n} \leq \psi(r) \stackrel{\text{def}}{=} \max\{\alpha, (r/C)^\beta\} .$$

It says that $|B(n, r)|$ is a growing function in $r$, and the growth rate is $\propto r^\beta$. We also allow an additive term $\alpha$ to cover the case that an $\alpha$ fraction of the stochastic gradients always need to be recalculated, regardless of the distance. We shall later illustrate why Assumption 2 holds in practice and why Assumption 2 holds in theory under reasonable data assumptions.

Our result of this section can be summarized as follows. Hiding $\|x_0 - x^*\|, L, C, \beta$ in the big-$O$ notion, and letting $T_{\text{time}}$ be the gradient complexity, we can modify GD so that it finds a point $x$ with

$$f(x) - f(x^*) \leq \begin{cases} O\big(\frac{\alpha}{T_{\text{time}}} + 2^{-\Omega(T_{\text{time}})^{1/3}}\big), & \text{if } \beta = 1; \\ \widetilde{O}\big(\frac{\alpha}{T_{\text{time}}} + T_{\text{time}}^{-1/(1-\beta)}\big), & \text{if } \beta \in (0, 1). \end{cases}$$

We emphasize that our modified algorithm does not need to know $\alpha$ or $\beta$.

---

[3]More generally, there is an "optimal" way to tweak the non-smooth objective to allow essentially any smooth-objective solver to apply, see [4].

**Algorithm 1** $\texttt{GD}^{\texttt{lin}}(f, x^{(0)}, S, C, D)$

---
**Input:** $f(x) = \frac{1}{n}\sum_{i=1}^n f_i(x)$ convex and $L$-smooth, starting vector $x^{(0)} \in \mathbb{R}^d$, number of epochs $S \geq 1$, parameters $C, D > 0$.
**Output:** vector $x \in \mathbb{R}^d$.
1: **for** $s \leftarrow 1$ **to** $S$ **do**
2:     $x_0 \leftarrow x^{(s-1)}$; $m \leftarrow \lceil (1 + \frac{C^2}{16D^2})^s \rceil$; and $\xi \leftarrow \frac{C}{m}$.
3:     $\mathbf{g} \leftarrow \vec{0}$ and $\mathbf{g}_i \leftarrow \vec{0}$ for each $i \in [n]$.
4:     **for** $k \leftarrow 0$ **to** $m - 1$ **do**
5:        Calculate $\Lambda_k \subseteq [n]$ from $x_0, \ldots, x_k$ according to Definition 3.1.
6:        **for** $i \in \Lambda_k$ **do**
7:           $\mathbf{g} \leftarrow \mathbf{g} + \frac{\nabla f_i(x_k) - \mathbf{g}_i}{n}$ and $\mathbf{g}_i \leftarrow \nabla f_i(x_k)$.
8:        $x_{k+1} \leftarrow x_k - \min\{\frac{\xi}{\|\mathbf{g}\|}, \frac{1}{L}\}\mathbf{g}$                           ⋄ *it satisfies* $\mathbf{g} = \nabla f(x_k)$
9:     $x^{(s)} \leftarrow x_m$;
10: **return** $x = x^{(S)}$.

---

## 3.1 Algorithm Description

In classical gradient descent (GD), starting from $x_0 \in \mathbb{R}^d$, one iteratively updates $x_{k+1} \leftarrow x_k - \frac{1}{L}\nabla f(x_k)$. We propose $\texttt{GD}^{\texttt{lin}}$ (see Algorithm 1) which, at a high level, differs from GD in two ways:

- It performs a truncated gradient descent with travel distance $\|x_k - x_{k+1}\| \leq \xi$ per step.
- It speeds up the process of calculating $\nabla f(x_k)$ by using the lingering of past gradients.

Formally, $\texttt{GD}^{\texttt{lin}}$ consists of $S$ epochs $s = 1, 2, \ldots, S$ of growing length $m = \lceil (1 + \frac{C^2}{16D^2})^s \rceil$. In each epoch, it starts with $x_0 \in \mathbb{R}^d$ and performs $m$ truncated gradient descent steps

$$x_{k+1} \leftarrow x_k - \min\{\tfrac{\xi}{\|\nabla f(x_k)\|}, \tfrac{1}{L}\} \cdot \nabla f(x_k) \ .$$

We choose $\xi = C/m$ to ensure that the worst-case travel distance $\|x_m - x_0\|$ is at most $m\xi = C$. (Recall that $r = C$ is the maximum distance so that $\psi(r) \leq 1$.)

In each iteration $k = 0, 1, \ldots, m-1$ of this epoch $s$, in order to calculate $\nabla f(x_k)$, $\texttt{GD}^{\texttt{lin}}$ constructs index sets $\Lambda_0, \Lambda_1, \ldots, \Lambda_{m-1} \subseteq [n]$ and recalculates only $\nabla f_i(x_k)$ for those $i \in \Lambda_k$. We formally introduce index sets below, and illustrate them in Figure 2(a).

**Definition 3.1.** *Given $x_0, x_1, \ldots, x_{m-1} \in \mathbb{R}^d$, we define index subsets $\Lambda_0, \ldots \Lambda_{m-1} \subseteq [n]$ as follows. Let $\Lambda_0 = [n]$. For each $k \in \{1, 2, \ldots, m-1\}$, if $(k_0, \ldots, k_t)$ is $k$'s lowbit sequence from Definition 3.2, then (recalling $k = k_t$)*

$$\Lambda_k \stackrel{\mathrm{def}}{=} \bigcup_{i=0}^{t-1} \big(B_{k_i}(k - k_i) \setminus B_{k_i}(k_{t-1} - k_i)\big)$$

*where*

$$B_k(r) \stackrel{\mathrm{def}}{=} \Lambda_k \cap B(x_k, r \cdot \xi) \ .$$

In the above definition, we have used the notion of "lowbit sequence" for a positive integer.[4]

**Definition 3.2.** *For positive integer $k$, let $\mathsf{lowbit}(k) \stackrel{\mathrm{def}}{=} 2^i$ where $i \geq 0$ is the maximum integer such that $k$ is integral multiple of $2^i$. For instance, $\mathsf{lowbit}(34) = 2$, $\mathsf{lowbit}(12) = 4$, and $\mathsf{lowbit}(8) = 8$.*

*Given positive integer $k$, let the **lowbit sequence** of $k$ be $(k_0, k_1, \ldots, k_t)$ where*

$$0 = k_0 < k_1 < \cdots < k_t = k \quad \text{and} \quad k_{i-1} = k_i - \mathsf{lowbit}(k_i) \ .$$

---
[4] If implemented in C++, we have $\mathsf{lowbit}(k) = $ `k & (-k)`.

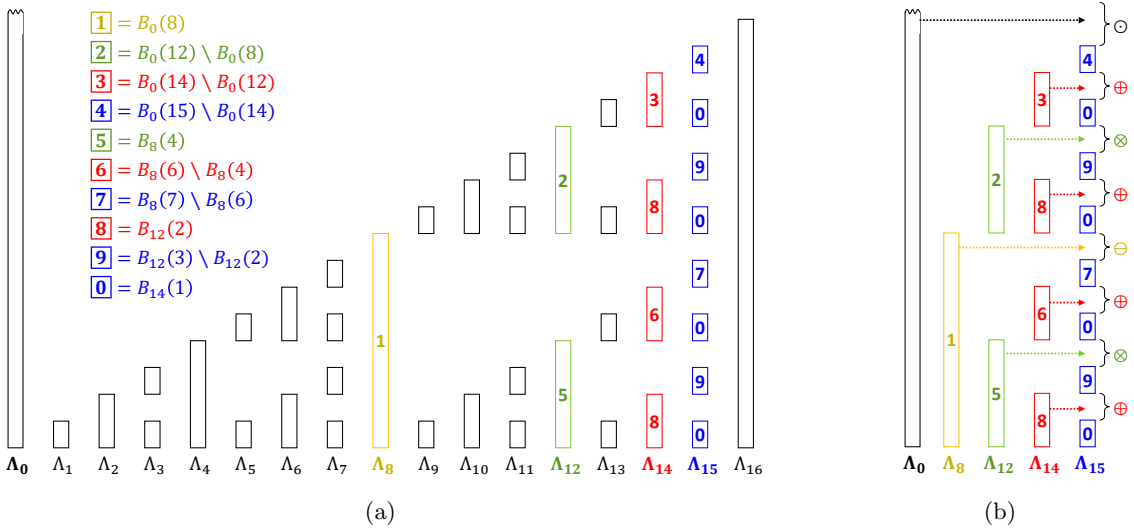

Figure 2: Illustration of index sets $\Lambda_k$

*For instance, the lowbit sequence of $45$ is $(0, 32, 40, 44, 45)$.*

## 3.2 Intuitions & Properties of Index Sets

We show in this paper that our construction of index sets satisfy the following three properties.

**Lemma 3.3.** *The construction of $\Lambda_0, \ldots, \Lambda_{m-1}$ ensures that $\mathbf{g} = \nabla f(x_k)$ in each iteration $k$.*

**Claim 3.4.** *The gradient complexity to construct $\Lambda_0, \ldots, \Lambda_{m-1}$ is $O\left(\frac{1}{n}\sum_{k=0}^{m-1}|\Lambda_k|\right)$ under Assumption 1. The space complexity is $O(n \log n)$.*

**Lemma 3.5.** *Under Assumption 2, we have $\frac{1}{n}\sum_{k=0}^{m-1}|\Lambda_k| \leq O(\alpha m + m^{1-\beta}\log^2 m)$ .*

At high level, Lemma 3.3 ensures that $\texttt{GD}^{\texttt{lin}}$ follows exactly the full gradient direction per iteration; Claim 3.4 and Lemma 3.5 together ensure that the total gradient complexity for this epoch is only $\widetilde{O}(m^{1-\beta}\log^2 m)$, as opposed to $O(m)$ if we always recalculate $\nabla f_1(x_k), \ldots, \nabla f_n(x_k)$.

Claim 3.4 is easy to verify. Indeed, for each $\Lambda_\ell$ that is calculated, we can sort its indices $j \in \Lambda_\ell$ in the increasing order of $\delta(x_k, j)$.[5] Now, whenever we calculate $B_{k_i}(k-k_i) \setminus B_{k_i}(k_{t-1}-k_i)$, we have already sorted the indices in $\Lambda_{k_i}$, so can directly retrieve those $j$ with $\delta(x_{k_i}, j) \in (k_{t-1}-k_i, k-k_i]$.

As for the space complexity, in any iteration $k$, we only need to store $\lceil \log_2 k \rceil$ index sets $\Lambda_\ell$ for $\ell < k$. For instance, when calculating $\Lambda_{15}$ (see Figure 2(b)), we only need to use $\Lambda_0, \Lambda_8, \Lambda_{12}, \Lambda_{14}$; and from $k = 16$ onwards, we no longer need to store $\Lambda_1, \ldots, \Lambda_{15}$.

Lemma 3.3 is technically involved to prove (see Appendix A.2), but we give a sketched proof by picture. Take $k = 15$ as an example. As illustrated by Figure 2(b), for every $j \in [n]$,

- If $j$ belongs to $\Lambda_{15}$—i.e., boxes $4, 0, 9, 7$ of Figure 2—

  We have calculated $\nabla f_j(x_k)$ so are fine.

---
[5]Calculating those lingering radii $\delta(x_k, j)$ require gradient complexity $|\Lambda_\ell|$ according to Assumption 1, and the time for sorting is negligible.

- If $j$ belongs to $\Lambda_{14} \setminus B_{14}(1)$—i.e., $\oplus$ region of Figure 2(b)—

  We have $\nabla f_j(x_{15}) = \nabla f_j(x_{14})$ because $\|x_{15} - x_{14}\| \leq \xi$ and $j \notin B_{14}(1)$. Therefore, we can safely retrieve $\mathbf{g}_j = \nabla f_j(x_{14})$ to represent $\nabla f_j(x_{15})$.

- If $j$ belongs to $\Lambda_{12} \setminus B_{12}(3)$—i.e., $\otimes$ region of Figure 2(b)—

  We have $\nabla f_j(x_{15}) = \nabla f_j(x_{12})$ for similar reason above. Also, the most recent update of $\mathbf{g}_j$ was at iteration 12, so we can safely retrieve $\mathbf{g}_j$ to represent $\nabla f_j(x_{15})$.

- And so on.

In sum, for all indices $j \in [n]$, we have $\mathbf{g}_j = \nabla f_j(x_k)$ so $\mathbf{g} = \frac{\mathbf{g}_1 + \cdots + \mathbf{g}_n}{n}$ equals $\nabla f(x_k)$.

Lemma 3.5 is also involved to prove (see Appendix A.3), but again should be intuitive from the picture. The indices in boxes $1, 2, 3, 4$ of Figure 2 are disjoint, and belong to $B(x_0, 15\xi)$, totaling at most $|B(x_0, 15\xi)| \leq n\psi(15\xi)$. The indices in boxes $5, 6, 7$ of Figure 2 are also disjoint, and belong to $B(x_8, 7\xi)$, totaling at most $|B(x_8, 7\xi)| \leq n\psi(7\xi)$. If we sum up the cardinality of these boxes by carefully grouping them in this manner, then we can prove Lemma 3.5 using Assumption 2.

### 3.3 Convergence Theorem

So far, Lemma 3.5 shows we can reduce the gradient complexity from $O(m)$ to $\widetilde{O}(m^{1-\beta})$ for every $m$ steps of gradient descent. Therefore, we wish to set $m$ as large as possible, or equivalently $\xi = C/m$ as small as possible. Unfortunately, when $\xi$ is too small, it will impact the performance of truncated gradient descent (see Lemma A.3 in appendix). This motivates us to start with a small value of $m$ and increase it epoch by epoch. Indeed, as the number of epoch grows, $f(x_0)$ becomes closer to the minimum $f(x^*)$, and thus we can choose smaller values of $\xi$.

Formally, we have (proved in Appendix A.5)

**Theorem 3.6.** *Given any $x^{(0)} \in \mathbb{R}^d$ and $D > 0$ that is an upper bound on $\|x^{(0)} - x^*\|$. Suppose Assumption 1 and 2 are satisfied with parameters $C \in (0, D], \alpha \in [0, 1], \beta \in (0, 1]$. Then, denoting by $m_s = \lceil (1 + \frac{C^2}{16D^2})^s \rceil$, we have that $\mathtt{GD}^{\mathtt{lin}}(f, x_0, S, C, D)$ outputs a point $x \in \mathbb{R}^d$ satisfying $f(x) - f(x^*) \leq \frac{4LD^2}{m_S}$ with gradient complexity $T_{\mathsf{time}} = O\big(\sum_{s=1}^S \alpha m_s + m_s^{1-\beta} \log^2 m_s\big)$.*

As simple corollaries, we have (proved in Appendix A.6)

**Theorem 3.7.** *In the setting of Theorem 3.6, given any $T \geq 1$, one can choose $S$ so that*

- *If $\beta = 1$, then $\mathtt{GD}^{\mathtt{lin}}$ finds a point $x$ in gradient complexity $T_{\mathsf{time}} = O(T)$ s.t.*

$$f(x) - f(x^*) \leq O\big(\tfrac{LD^4}{C^2} \cdot \tfrac{\alpha}{T}\big) + \tfrac{LD^2}{2^{\Omega(C^2T/D^2)^{1/3}}} \ .$$

- *If $\beta \in (0,1)$ is constant, $\mathtt{GD}^{\mathtt{lin}}$ finds a point $x$ in gradient complexity $T_{\mathsf{time}} = O(T \log^2 T)$ s.t.*

$$f(x) - f(x^*) \leq O\big(\tfrac{LD^4}{C^2} \cdot \tfrac{\alpha}{T} + \tfrac{LD^{2+\frac{2}{1-\beta}}}{C^{\frac{2}{1-\beta}}} \cdot \tfrac{1}{T^{1/(1-\beta)}}\big) \ .$$

We remark here if $\psi(r) = 1$ (so there is no lingering effect for gradients), we can choose $C = D$ and $\beta = 1$; in this case $\mathtt{GD}^{\mathtt{lin}}$ gives back the convergence $f(x) - f(x^*) \leq O\big(\frac{LD^2}{T}\big)$ of GD.

## 4 SVRG with Lingering Radius

In this section, we use Assumption 1 to improve the running time of SVRG [21, 38], one of the most widely applied stochastic gradient methods in large-scale settings. The purpose of this section is to construct an algorithm that works well *in practice*: to (1) work for any possible lingering radii $\delta(x, i)$, (2) be identical to SVRG if $\delta(x, i) \equiv 0$, and (3) be faster than SVRG when $\delta(x, i)$ is large.

Recall how the SVRG method works. Each *epoch* of SVRG consists of $m$ iterations ($m = 2n$ in practice). Each epoch starts with a point $x_0$ (known as the *snapshot*) where the full gradient $\nabla f(x_0)$ is computed exactly. In each iteration $k = 0, 1, \ldots, m-1$ of this epoch, SVRG updates $x_{k+1} \leftarrow x_k - \eta \mathbf{g}$ where $\eta > 0$ is the learning rate and $\mathbf{g}$ is the gradient estimator $\mathbf{g} = \nabla f(x_0) + \nabla f_i(x_k) - \nabla f_i(x_0)$ for some $i$ randomly drawn from $[n]$. Note that it satisfies $\mathbb{E}_i[\mathbf{g}] = \nabla f(x_k)$ so $\mathbf{g}$ is an unbiased estimator of the gradient. In the next epoch, SVRG starts with $x_m$ of the previous epoch.[6] We denote by $x^{(s)}$ the value of $x_0$ at the beginning of epoch $s = 0, 1, 2, \ldots, S-1$.

## 4.1 Algorithm Description

Our algorithm `SVRG`[lin] maintains *disjoint* subsets $H_s \subseteq [n]$, where each $H_s$ includes the set of the indices $i$ whose gradients $\nabla f_i(x^{(s)})$ from epoch $s$ can still be safely reused at present.

At the starting point $x_0$ of an epoch $s$, we let $H_s = [n] \setminus (H_0 \cup \cdots \cup H_{s-1})$ and re-calculate gradients $\nabla f_i(x_0)$ only for $i \in H_s$; the remaining ones can be loaded from the memory. This computes the full gradient $\nabla f(x_0)$. Then, we denote by $m = 2|H_s|$ and perform only $m$ iterations within epoch $s$. We next discuss how to perform update $x_k \to x_{k+1}$ and maintain $\{H_s\}_s$ during each iteration.

- In each iteration $k$ of this epoch, we claim that $\nabla f_i(x_k) = \nabla f_i(x_0)$ for every $i \in H_0 \cup \cdots \cup H_s$.[7] Thus, we can uniformly sample $i$ from $[n] \setminus (H_0 \cup \cdots \cup H_s)$, and construct an unbiased estimator

$$\mathbf{g} \leftarrow \nabla f(x_0) + \left(1 - \frac{\sum_{s'=0}^{s} |H_{s'}|}{n}\right) [\nabla f_i(x_k) - \nabla f_i(x_0)]$$

of the true gradient $\nabla f(x_k)$. Then, we update $x_{k+1} \leftarrow x_k - \eta \mathbf{g}$ the same way as SVRG. We emphasize that the above choice of $\mathbf{g}$ reduces its variance (because there are fewer random choices), and it is known that reducing variance leads to faster running time [21].

- As for how to maintain $\{H_s\}_s$, in each iteration $k$ after $x_{k+1}$ is computed, for every $s' \leq s$, we wish to remove those indices $i \in H_{s'}$ such that the current position $x$ lies outside of the lingering radius of $i$, i.e., $\delta(x^{(s)}, i) < \|x - x^{(s)}\|$. To efficiently implement this, we need to make sure that whenever $H_{s'}$ is constructed (at the beginning of epoch $s'$), the algorithm sorts all the indices $i \in H_{s'}$ by increasing order of $\delta(x^{(s')}, i)$. We include implementation details in Appendix B.1.

## 4.2 SCSG with Lingering Radius

When $n$ is extremely large, it can be expensive to compute full gradient at snapshots, so a variant of SVRG is sometimes applied in practice. That is, at each snapshot $x_0$, instead of calculating $\nabla f(x_0)$, one can approximate it by a batch average $\frac{1}{|S|} \sum_{i \in S} \nabla f_i(x_0)$ for a sufficiently large random subset $S$ of $[n]$. Then, the length of an epoch is also changed from $m = 2n$ to $m = 2|S|$. This method is studied by [17, 22, 23], and we refer to it as SCSG due to [22, 23].

Our algorithm `SVRG`[lin] can be easily extended to this setting, with the following modifications:

---

[6]Some authors use the average of $x_1, \ldots, x_m$ to start the next epoch, but we choose this simpler version.

[7]This is because for every $i \in H_s$, by definition of $H_s$ we have $\nabla f_i(x_k) = \nabla f_i(x^{(s)}) = \nabla f_i(x_0)$; for every $i \in H_{s'}$ where $s' < s$, we know $\nabla f_i(x_k) = \nabla f_i(x^{(s')})$ but we also have $\nabla f_i(x_0) = \nabla f_i(x^{(s')})$ (because otherwise $i$ would have been removed from $H_{s'}$).

**Algorithm 2** $\text{SVRG}^{\text{lin}}(f, x^{(0)}, \eta, S)$

**Input:** $f(x) = \frac{1}{n} \sum_{i=1}^{n} f_i(x)$, vector $x^{(0)} \in \mathbb{R}^d$, learning rate $\eta > 0$, number of epochs $S \geq 1$.
**Output:** vector $x \in \mathbb{R}^d$.

1: **for** $s \leftarrow 0$ **to** $S - 1$ **do**
2: $\quad x_0 \leftarrow x^{(s)}; \quad H_s \leftarrow [n] \setminus (H_0 \cup \cdots \cup H_{s-1}); \quad$ and $m \leftarrow 2|H_s|$.
3: $\quad$ compute full gradient $\nabla f(x_0)$ according to
$$\nabla f(x_0) = \tfrac{1}{n} \left[ \sum_{s'=0}^{s-1} \sum_{i \in H_{s'}} \nabla f_i(x^{(s')}) + \sum_{i \in H_s} \nabla f_i(x_0) \right].$$
4: $\quad$ **for** $k \leftarrow 0$ **to** $m - 1$ **do**
5: $\quad\quad$ **if** $H_0 \cup \cdots \cup H_s = [n]$ **then**
6: $\quad\quad\quad \mathbf{g} \leftarrow \nabla f(x_0)$.
7: $\quad\quad$ **else**
8: $\quad\quad\quad$ randomly draw $i \in [n] \setminus (H_0 \cup \cdots \cup H_s)$.
9: $\quad\quad\quad \mathbf{g} \leftarrow \nabla f(x_0) + \left( 1 - \frac{\sum_{s'=0}^{s} |H_{s'}|}{n} \right) [\nabla f_i(x_k) - \nabla f_i(x_0)]$.
10: $\quad\quad x_{k+1} \leftarrow x_k - \eta \mathbf{g}$.
11: $\quad\quad$ **for all** $s' \leq s$ and $i \in H_{s'}$ such that $\delta(x^{(s')}, i) < \|x^{(s')} - x_{k+1}\|$ **do**
12: $\quad\quad\quad H_{s'} \leftarrow H_{s'} \setminus \{i\}$.
13: $\quad x^{(s+1)} \leftarrow x_m$.
14: **return** $x = x^{(S)}$.

- We define parameter $\overline{m}_s = \min\{n, \overline{m}_0 \cdot 2^s\}$, where $\overline{m}_0$ is a given input (allegedly the length of the first epoch).

- We replace Line 2 of $\text{SVRG}^{\text{lin}}$ as follows. If there are more than $\overline{m}_s$ elements in $[n] \setminus (H_0 \cup \cdots \cup H_{s-1})$, then

    $H_s \leftarrow$ a random sample of $\overline{m}_s$ indices from $[n] \setminus (H_0 \cup \cdots \cup H_{s-1})$.

    Otherwise, set $H_s$ in the same way as in $\text{SVRG}^{\text{lin}}$.

- We replace Line 3 of $\text{SVRG}^{\text{lin}}$, the full gradient computation, by an estimate
$$\nabla f(x_0) \approx \frac{1}{n} \left[ \sum_{s'=0}^{s-1} \sum_{i \in H_{s'}} \nabla f_i(x^{(s')}) + \frac{|[n] \setminus (H_0 \cup \cdots \cup H_{s-1})|}{|H_s|} \sum_{i \in H_s} \nabla f_i(x_0) \right].$$

    It can be computed using $|H_s| \leq \overline{m}_s$ computations of new gradients.

We call this algorithm $\text{SCSG}^{\text{lin}}$ and also report its practical performance in our experiments. We note that having epoch size to grow exponentially was recommended for instance by the authors of SCSG [23] and others [5, 25].

## 5 Experiments on Packing LP

In this section, we construct a revenue maximization LP (2.1) using the publicly accessible dataset of Yahoo! Front Page Today Module [7, 24]. Based on this real-life dataset, we validate Assumption 2 and our motivation behind lingering gradients. We also test the performance of $\text{SVRG}^{\text{lin}}$ from Section 4 and $\text{SCSG}^{\text{lin}}$ from Section 4.2 on optimizing this LP.

## 5.1 Experiment Setup

We use part of the Today Module dataset corresponding to May 1, 2009. There are $d = 50$ articles, which we view as resources, and $n \approx 4.6$ million users. We estimate $p_{i,j}$ following the hybrid model in [24]. While Li et al. [24] consider the online recommendation problem without any constraints on the total traffic that each article receives, we consider the *offline* LP problem (2.1) with resource capacity constraints. In practice, recommendation systems with resource constraints can better control the public exposure of any ads or recommendations [39].

In addition to estimating $p_{i,j}$ from data, we generate other synthetic parameters in order to make the LP problem (2.1) non-trivial to solve. From a high level, we want (i) some resources to have positive remaining capacities under optimal LP solutions, so that the LP is feasible (when (2.1) is infeasible due to the equality constraints, the revenue-maximization problem becomes trivial because we can sell all the inventories); (ii) some resources to have zero remaining capacities under optimal LP solutions, so that the optimal dual solution is not a (trivial) zero vector. Specifically,

- We arbitrarily pick a resource $k \in [d]$, and assign it infinity capacity $b_k > n$ with relatively small revenue value $r_k = 0.05$.
- For other resources $i \in [d]$, we randomly draw $r_i$ from a uniform distribution over $[0.05, 0.95]$, and set $b_i = 0.01n/d$.
- We choose $\mu = 10^{-5}$ as the regularization error.
- For each algorithm, we tune learning rates from the set $\eta \in \{10^{-k}, 3 \times 10^{-k}, 5 \times 10^{-k}\}$, and report the best-tuned performance.

Finally, we note that the dual objective (2.2) is constrained optimization with $x \geq 0$. Although we specified our algorithm SVRG$^{\text{lin}}$ (for notational simplicity) without constraints on $x$, it is a simple exercise to generalize it (as well as classical methods SVRG, SAGA) into the constrained setting. Namely, in each step, if the new point $x_{k+1}$ moves out of the constraint, then project it to the closest point on the constraint. This is known as the proximal setting of first-order method, and see for instance the analysis of proximal SVRG of [37].

We discuss implementation details of SVRG$^{\text{lin}}$ and SCSG$^{\text{lin}}$ in Appendix B.

## 5.2 Illustration of Lingering Radius

We calculate lingering radii on the dual problem (2.2). Let $\theta > 0$ be a parameter large enough so that $e^{-\theta}$ can be viewed as zero. (For instance, $\theta = 20$ gives $e^{-20} \approx 2 \times 10^{-9}$.) Then, for each point $x \in \mathbb{R}_{\geq 0}$ and index $i \in [n]$, we let

$$\delta(x, i) = \min_{j \in [n], j \neq j^*} \frac{(r_{j^*} - x_{j^*})p_{i,j^*} - (r_j - x_j)p_{i,j} - \theta \overline{p_i}\mu}{p_{i,j^*} + p_{i,j}}$$

where

$$j^* = \arg\max_{j \in [n]} \{(r_j - x_j)p_{i,j}\}.$$

It is now a simple exercise to verify that, denoting by $\mathbf{e}_j$ the $j$-th basis unit vector, then[8]

$$\nabla f_i(y) \approx (b_1, \ldots, b_d) + np_{i,j^*}\mathbf{e}_{j^*} \quad \text{for every} \quad \|y - x\|_\infty \leq \delta(x, i) \ .$$

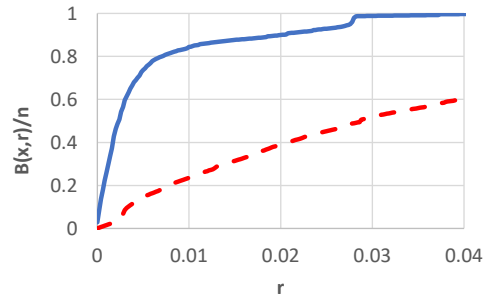

Figure 3: $|B(x, r)|/n$ as a function of $r$ for packing LP. We choose $\theta = 5$. Dashed curve is when $x = \vec{0}$, and solid curve is when $x$ is near the optimum.

---
[8]For any other coordinate $j \neq j^*$, it satisfies $\frac{e^{(r_j - y_j)p_{i,j}/(\overline{p_i}\mu)}}{e^{(r_{j^*} - y_{j^*})p_{i,j^*}/(\overline{p_i}\mu)}} \leq e^{-\theta}$ and hence is negligible.

In Figure 3, we plot $|B(x,r)| = |\{j \in [n] \,|\, \delta(x,j) < r\}|$ as an increasing function of $r$. We see that for practical data, $|B(x,r)|/n$ is indeed bounded above by some increasing function $\psi(\cdot)$.

We justify Figure 3 as follows. For any point $x$ and customer $i$, recall from (2.3) that $y_{i,j} \propto \exp(\frac{(r_j-x_j)p_{i,j}}{\bar{p}_i \mu})$ approximately captures the index $j = j^*$ which maximizes the exponent. If $\mu$ is small (recall small $\mu$ gives more accurate solutions to primal LP (2.1)), then for most customers, $(y_{i,1}, \ldots, y_{i,d})$ is approximately a unit vector $\mathbf{e}_{j^*}$, meaning we assign customer $i$ to resource $j$ with high probability. Now, as long as $x$ stays in the lingering radius of $i$, we still offer customer $i$ the same resource $j^*$. Naturally, when customers are sufficiently diverse— which is usually the case in practice— one should expect the lingering radii to be evenly distributed, and thus $|B(x,r)|$ can behave like Figure 3.

*Remark* 5.1. This $\delta(x,i)$ differs from our definition in Section 2 in two ways. First, it ensures $\nabla f_i(y) \approx \nabla f_i(x)$ as opposed to exact equality; for practical purposes this is no big issue, and we choose $\theta = 5$ in our experiments. Second, $\|y - x\|_\infty \leq \delta(x,i)$ gives a bigger "safe region" than $\|y - x\| \leq \delta(x,i)$; thus, when implementing `SVRG`<sup>lin</sup>, we adopt $\|\cdot\|_\infty$ as the norm of choice.

*Remark* 5.2. Figure 3 also confirms that Assumption 2 holds in practice. Recall Assumption 2 was used in proving the theoretical performance of `GD`<sup>lin</sup>. We see that indeed $|B(x,r)|$ grows increasingly in $r$.

## 5.3 Performance Comparison

We consider solving the dual problem (2.2). In Figure 4(a) and 4(b), we plot the optimization error of (2.2) as a function $\#\mathrm{grad}/n$, the number of stochastic gradient computations divided by $n$, also known as $\#passes\ of\ dataset$.

Figure 4(a) compares our `SVRG`<sup>lin</sup> to `SVRG` and `SAGA` (each for 3 best tuned learning rates), and Figure 4(b) compares our `SCSG`<sup>lin</sup> to `SCSG` (also with 3 best tuned learning rates).[9] We can see `SVRG`<sup>lin</sup> is close to `SVRG` or `SAGA` during the first 5-7 passes of the data. This is because initially, $x$ moves fast and cannot usually stay in the lingering radii for most indices $i$. After that period, `SVRG`<sup>lin</sup> requires a dramatically smaller number of gradient computations, as $x$ moves slower and slower, becoming more easily to stay in the lingering radii. It is interesting to note that `SVRG`<sup>lin</sup> does not significantly improve the optimization error as a function of number of epochs; the improvement primarily lies in improving the number of gradient computations per epoch. The comparison is similar for `SCSG`<sup>lin</sup> vs. `SCSG`.

## 5.4 Performance Comparison on Primal LP Objective

When $f(x)$ is only approximately minimum, the optimization error of the dual (2.2) does not represent the error for the primal LP (2.1). Therefore, the more interesting notion is the *primal error*, defined as

$$\left[\mathtt{OPT} - \sum_{j \in [d]} r_j \min(b_j, \sum_{i \in [n]} p_{i,j} y_{i,j})\right] / \mathtt{OPT},$$

where `OPT` is the optimal objective value of (2.1), and $y$ is given by (2.3). This primal error captures the inefficiency caused by the in-feasibility of $y$. Indeed, when $x$ is not the exact minimizer, the amount of demand assigned to a resource $j$ may exceed its capacity $b_j$. If this happens, in the above expression, we measure the primal objective with respect to $y$ by truncating all the demand that exceeds its capacity.

---

[9] Each epoch of `SVRG` consists of a full gradient computation and $2n$ iterations, totaling $3n$ computations of (new) stochastic gradients. (We do not count the computation of $\nabla f_i(0)$ at $x = 0$.) Each epoch of `SCSG` needs to compute a batch gradient of size $|S|$, followed by $2|S|$ iterations each computing 2 new stochastic gradients.

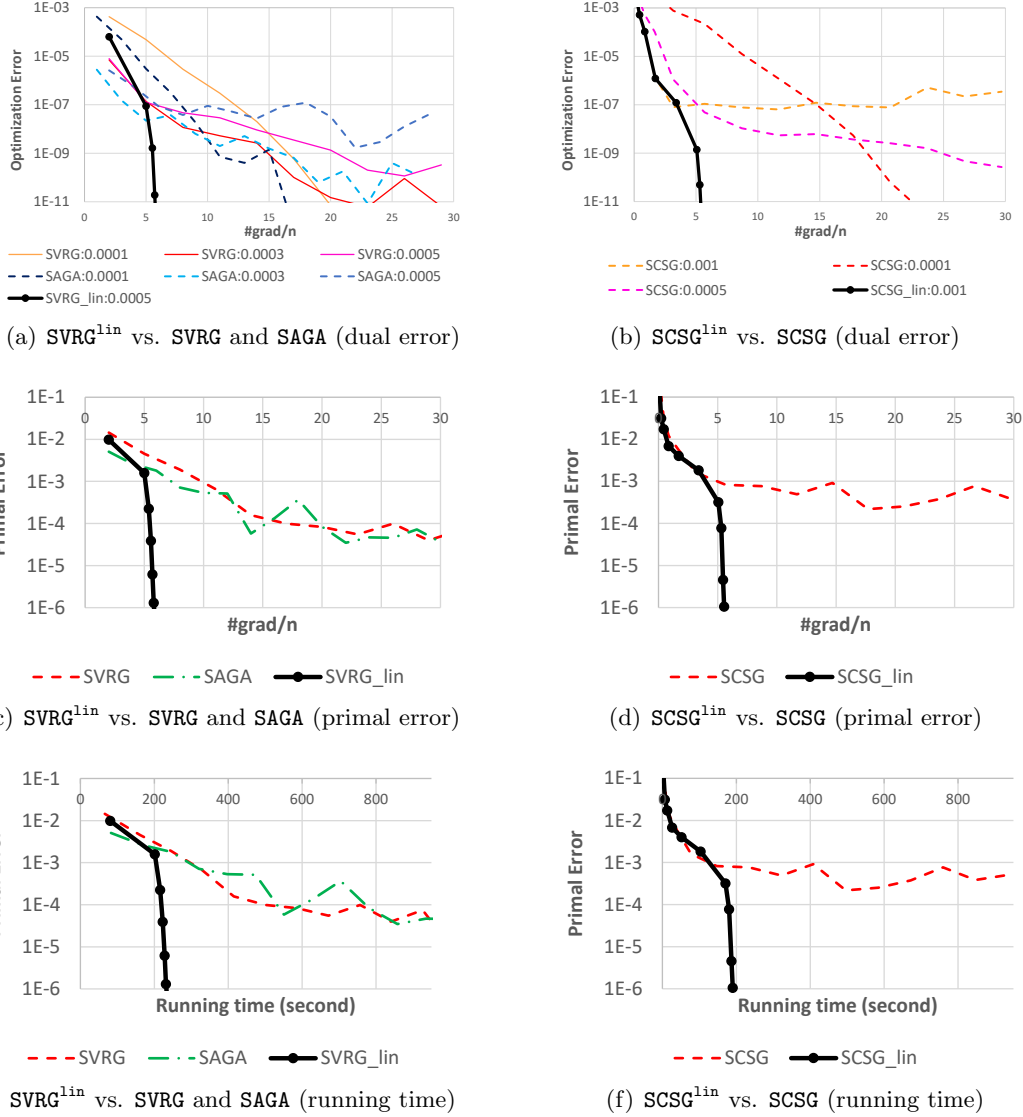

Figure 4: Performance comparison for revenue management LP.

Figure 4(c) compares $\texttt{SVRG}^{\texttt{lin}}$ to $\texttt{SVRG}$ and $\texttt{SAGA}$ in the primal error, while Figure 4(d) compares $\texttt{SCSG}^{\texttt{lin}}$ to $\texttt{SCSG}$. (Here for simplicity, we have only plotted those prior works with respect to the best tuned learning rate.) We can see that standard stochastic descent algorithms need to spend more than 30 passes of data in order to achieve less than $10^{-4}$ primal error, while $\texttt{SVRG}^{\texttt{lin}}$ and $\texttt{SCSG}^{\texttt{lin}}$ converge to $10^{-6}$ within no more than 6 passes of the data. Note that the primal error also contains the loss caused by regularization (which recall we have chosen $\mu = 10^{-5}$); our algorithms quickly achieve small primal errors comparable to $\mu$.

In Figure 4(e) and 4(f), we also compare the running time of the algorithms.

# 6 Experiments on SVM

In this section, we construct an SVM objective for the binary classification task on the real-life libsvm Adult dataset [12].

## 6.1 Experiment Setup

In the Adult dataset there are $n = 32,561$ datapoints and $d = 123$ dimensions. We re-scale the data points by a common constant so that their average Euclidean norm is 1. We choose $\lambda = 1/n$ as the regularizer weight for the SVM objective.

- We have implemented the vanilla PEGASOS method of Shalev-Shwartz et al. [33] which does not require any parameter tuning and directly applies to the original SVM objective (2.4).

- For existing methods GD, SVRG, SCSG and our new methods $\texttt{GD}^{\texttt{lin}}, \texttt{SVRG}^{\texttt{lin}}, \texttt{SCSG}^{\texttt{lin}}$, we run them both on the original SVM objective (2.4) (denoted by $\mu = 0$) as well as the smooth SVM objective (2.5) (using $\mu = 0.01$). We point out that the theory for many of these methods do not apply to the original non-smooth SVM objective, but in practice this is not an issue.

- For each algorithm (except PEGASOS), we tune learning rates from the set $\eta \in \{10^{-k}, 2.5 \times 10^{-k}, 5 \times 10^{-k}, 7.5 \times 10^{-k}\}$, and report the best-tuned performance.

- When plotting the training performance of any of these methods, we only show the vanilla SVM objective (2.4), but never the smooth alternative.[10]

We discuss implementation details of $\texttt{SVRG}^{\texttt{lin}}$, $\texttt{SCSG}^{\texttt{lin}}$, and $\texttt{GD}^{\texttt{lin}}$ in Appendix B.

## 6.2 Illustration of Lingering Radius

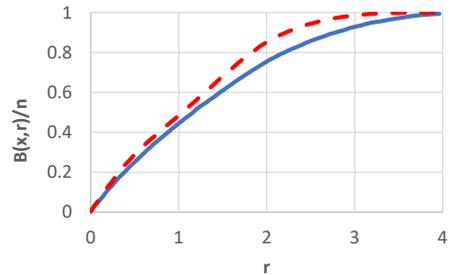

Figure 5: $|B(x,r)|/n$ as a function of $r$ for SVM. Dashed curve is when $x = \vec{0}$, and solid curve is when $x$ is near the optimum.

The calculation of lingering radii is trivial for SVM. Recall for the smoothed SVM objective, each stochastic gradient (ignoring the regularizer term $\frac{\lambda}{2}\|x\|^2$) can be explicitly written as[11]

$$\nabla f_i^\mu(x) = \begin{cases} 0, & b_i\langle x, a_i\rangle \geq 1; \\ -b_i a_i, & b_i\langle x, a_i\rangle \leq 1 - \mu; \\ \frac{-b_i a_i}{\mu}(1 - b_i\langle x, a_i\rangle), & \text{otherwise.} \end{cases} \quad (6.1)$$

From this definition, we can calculate $\delta(x,i)$ as follows. If $b_i\langle x, a_i\rangle \geq 1$ then we set $\delta(x,i) = \frac{b_i\langle x, a_i\rangle - 1}{\|a_i\|}$; if $b_i\langle x, a_i\rangle \leq 1 - \mu$ then we set $\delta(x,i) = \frac{1 - \mu - b_i\langle x, a_i\rangle}{\|a_i\|}$; and if otherwise, we set $\delta(x,i) = +\infty$. It is clear from the gradient formula that this definition is valid.

In Figure 5, we plot $|B(x,r)|/n$ as a function of $r$. (We choose the smooth version with $\mu = 0.01$ and the non-smooth version is only better.) We see that for practical data, $|B(x,r)|/n$ is indeed bounded from above by some almost-linear function in $r$ as theory predicts (see Section 7).

---

[10]In machine learning, the smooth SVM objective (2.5) is usually viewed as an auxiliary objective whose purpose is to help minimize the original SVM objective (2.4).

[11]The gradient of the regularizer is $\lambda x$ so can be calculated efficiently without the necessity of reading the dataset. For such reason, we do not need to consider it while calculating the lingering radius.

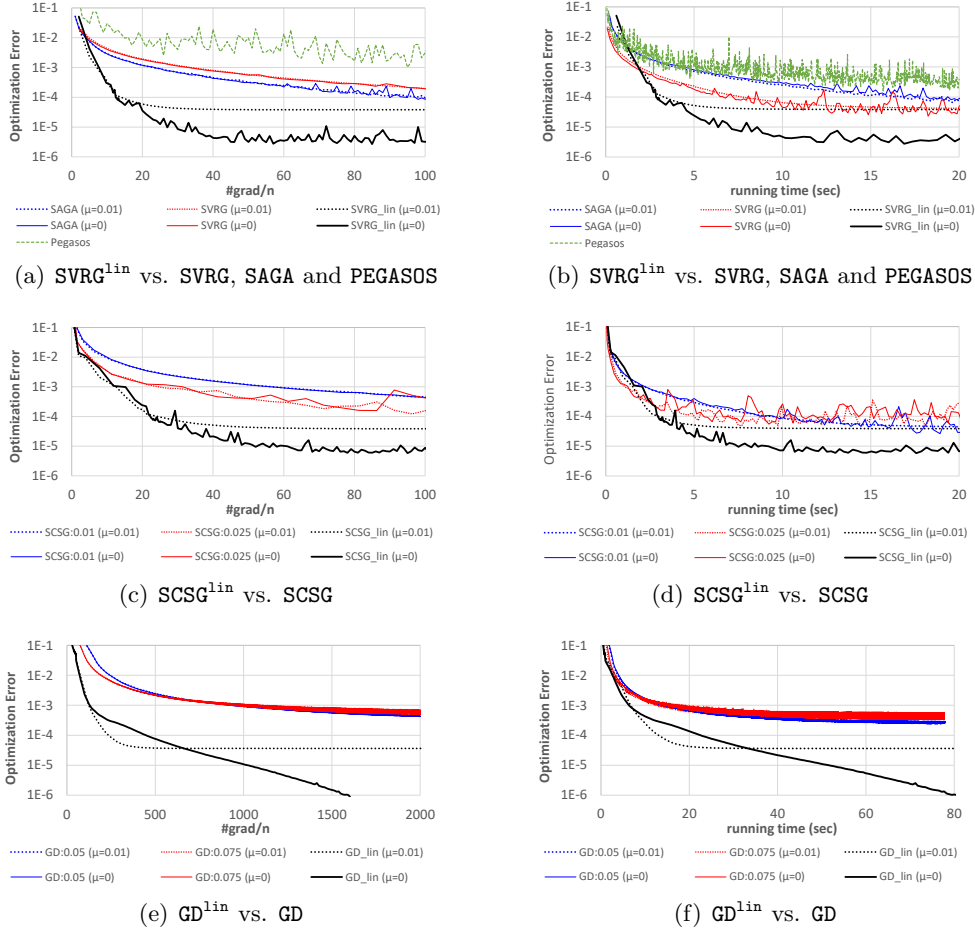

Figure 6: Performance comparison for training SVM.

## 6.3 Performance Comparison

We consider two types of performance charts. The first is the optimization error of the original SVM objective (2.4) as a function #grad/$n$, also known as #*passes of dataset*; the second is the optimization error as a function of the running time.

- Figure 6(a) and 6(b) compare SVRG$^{lin}$ to SVRG and SAGA (each with best tuned learning rates) as well as to PEGASOS. Due to the non-smooth nature of SVM, existing methods all have a hard time achieving high accuracy such as $10^{-5}$ error; in contrast, SVRG$^{lin}$ can easily achieve this accuracy within 30 passes of the data). Note also, if an algorithm is trained on the smoothed SVM objective (2.5) with $\mu > 0$, then its optimization error does not converge to zero on the original SVM objective (2.4).

- Figure 6(c) and 6(d) compare our SCSG$^{lin}$ (with best tuned learning rate) to SCSG (with two of the best tuned learning rates). Again, there is a clear performance gain to take into account the lingering of gradients. However, for this task of SVM, SCSG$^{lin}$ and SCSG do not seem to outperform the corresponding SVRG$^{lin}$ and SVRG.

- Figure 6(e) and 6(f) compare our GD$^{lin}$ (with best tuned learning rate) to GD (with two of the best tuned learning rates). Recall that GD$^{lin}$ serves as a *theoretical* evidence that exploiting the lingering of gradients can drastically improve the *convergence rate* of the method. From these two plots, one can verify that it is indeed so, and the performance gain is very similar to what theory predicts in Figure 1.

# 7 Theoretical Evidence for Assumption 2

Recall when proving our theoretical result in Section 3, we have made Assumption 2 which says $|B(x,r)|$ is bounded by an explicit increasing function in $r$. In this section, in the context of SVM, we demonstrate why it is so under some natural randomness assumption of the data.

**Assumption 3.** *Suppose there exists parameters $\sigma > 0$ and $\kappa \geq 1$ such that the following holds. Each data point $a_i \in \mathbb{R}^d$ is drawn independently from a Gaussian distribution $a_i \sim \mathcal{N}(\mu_i, \Sigma_i)$ with the property that $\Sigma_i \preceq \frac{\sigma^2}{d}\mathbf{I}$ and $\sup_{x \in \mathbb{R}^d} \left\{ \frac{|x^\top \mu_i|}{\sqrt{x^\top \Sigma_i x}} \right\} \leq \kappa$.*

As a simple example, if $\Sigma_i = \mathbf{I}$ and $\|\mu_i\| \leq 10$, then one can choose $\sigma = \sqrt{d}$ and $\kappa = 10$.

*Remark* 7.1. Assumption 3 is quite natural in the following sense. The bound $\Sigma_i \preceq \frac{\sigma^2}{d}\mathbf{I}$ ensures that the data point $a_i$ has a bounded Euclidean norm. The bound $\frac{|x^\top \mu_i|}{\sqrt{x^\top \Sigma_i x}} \leq \kappa$ ensures that $a_i$ cannot be very adversarial: the projection of $\mu$ to any direction $x$ must be bounded by the amount of randomness in that direction.

We main theorem of this section shows that Assumption 2 holds with $\beta = 1$ when data is sufficiently random:

**Theorem 7.2.** *Under Assumption 3, for every $\alpha_0 \in (0, 1)$ and $B > 0$, as long as $n \geq \Omega(\frac{d}{\alpha_0} \log \frac{B\sigma\kappa}{\alpha_0})$, with probability at least $1 - e^{-\Omega(\alpha_0 n)}$ it satisfies that for every $x \in \mathbb{R}^d$ with $\|x\| \leq B$ and $r \geq 0$,*

$$\frac{|B(x,r)|}{n} \leq \max\left\{ O(\sigma\kappa \log \frac{1}{\alpha_0}) \cdot r, \quad O\big(\mu\kappa \log \frac{1}{\alpha_0} + \alpha_0\big) \right\} .$$

## 7.1 Proof of Theorem 7.2

The proof of Theorem 7.2 consists of several steps of careful probabilistic derivations. In the first lemma below, we fix some vector $x \in \mathbb{R}^d$, fix some radius $r \geq 0$, and fix some data point $i \in [n]$. Using the randomness of $a_i$, this lemma upper bounds the probability for the value $b_i \langle x, a_i \rangle$ to be in a "bad" interval and $\|a_i\|$ being too large.

**Lemma 7.3.** *Under Assumption 3, for every $\alpha_0 \in (0, 1)$ there exists parameter $\xi \geq 1$ such that, for every $i \in [n]$, every $x \in \mathbb{R}^d$ and every $r \geq 0$, it satisfies*

$$\Pr\big[b_i\langle x, a_i\rangle \in [1 - 2r\xi\sigma - \mu, 1 + 2r\xi\sigma]\big] + \Pr[\|a_i\| \geq \xi\sigma] \leq O\big((r\sigma + \mu)\kappa \log \frac{1}{\alpha_0}\big) + \frac{\alpha_0}{5} .$$

*Proof of Lemma 7.3.* We first note by tail bound for chi-square distribution that $\Pr[\|a_i\| > \xi\sigma] \leq e^{-\xi^2 d/5}$. As for the other probability,

$$\Pr\big[b_i\langle x, a_i\rangle \in [1 - 2r\xi\sigma - \mu, 1 + 2r\xi\sigma]\big]$$
$$\leq \Pr\big[\langle x, a_i\rangle \in [1 - 2r\xi\sigma - \mu, 1 + 2r\xi\sigma]\big] + \Pr\big[-\langle x, a_i\rangle \in [1 - 2r\xi\sigma - \mu, 1 + 2r\xi\sigma]\big]$$

Now, notice that $\langle x, a_i \rangle \sim \mathcal{N}(x^\top \mu, x^\top \Sigma x)$. Therefore, define $g = \frac{\langle x, a_i\rangle - x^\top \mu}{\sqrt{x^\top \Sigma x}}$, it follows from distribution $g \sim \mathcal{N}(0,1)$. Therefore

$$\Xi \stackrel{\text{def}}{=} \Pr\big[\langle x, a_i\rangle \in [1 - 2r\xi\sigma - \mu, 1 + 2r\xi\sigma]\big] = \Pr\left[g \in \Big[\frac{1 - x^\top \mu}{\sqrt{x^\top \Sigma x}} - \frac{2r\xi\sigma + \mu}{\sqrt{x^\top \Sigma x}}, \frac{1 - x^\top \mu}{\sqrt{x^\top \Sigma x}} + \frac{2r\xi\sigma}{\sqrt{x^\top \Sigma x}}\Big]\right] .$$

On one hand, by property of standard Gaussian, there exists $t$ such that the above probability

$$\Xi = \Pr\left[g \in \left[t, t + \frac{4r\xi\sigma + \mu}{\sqrt{x^\top \Sigma x}}\right]\right] \leq \frac{1}{\sqrt{2\pi}} \frac{4r\xi\sigma + \mu}{\sqrt{x^\top \Sigma x}} \ .$$

On the other hand, by our assumption $\frac{|x^\top \mu|}{\sqrt{x^\top \Sigma x}} \leq \kappa$. If $\frac{1}{\sqrt{x^\top \Sigma x}} \geq 2\kappa$, then we know that $g$ in the above range implies $g \geq \frac{1/2}{\sqrt{x^\top \Sigma x}} - \frac{2r\xi\sigma + \mu}{\sqrt{x^\top \Sigma x}}$. Thus, if both $\frac{1}{\sqrt{x^\top \Sigma x}} \geq 2\kappa$ and $2r\xi\sigma + \mu \leq 1/4$ are satisfied, the above probability is at most

$$\Xi \leq \Pr\left[g \geq \frac{1 - x^\top \mu}{\sqrt{x^\top \Sigma x}} - \frac{2r\xi\sigma + \mu}{\sqrt{x^\top \Sigma x}}\right] \leq \Pr\left[g \geq \frac{1/2}{\sqrt{x^\top \Sigma x}} - \frac{2r\xi\sigma + \mu}{\sqrt{x^\top \Sigma x}}\right] \leq \Pr\left[g \geq \frac{1/4}{\sqrt{x^\top \Sigma x}}\right] \leq e^{-1/(32x^\top \Sigma x)} \ .$$

Finally, let $\nu \geq 2\kappa$ be a parameter to be chosen later. Consider the following two possibilities: $\frac{1}{\sqrt{x^\top \Sigma x}} \geq \nu$ and $\frac{1}{\sqrt{x^\top \Sigma x}} \leq \nu$. In the former case we have $\Xi \leq e^{-\Omega(\nu^2)}$; and in the latter case we have $\Xi \leq O((r\xi\sigma + \mu)\nu)$. In sum, we have

$$\Pr\left[\langle x, a_i \rangle \in [1 - 2r\xi\sigma - \mu, 1 + 2r\xi\sigma]\right] \leq O((r\xi\sigma + \mu)\nu) + e^{-\Omega(\nu^2)} \ .$$

For similar reason, we also have

$$\Pr\left[-\langle x, a_i \rangle \in [1 - r\xi\sigma - \mu, 1 + r\xi\sigma]\right] \leq O((r\xi\sigma + \mu)\nu) + e^{-\Omega(\nu^2)} \ .$$

Together, we conclude that for any parameter $\xi \geq 1$ and $\nu \geq 2\kappa$,

$$\Pr\left[b_i \langle x, a_i \rangle \in [1 - 2r\xi\sigma - \mu, 1 + 2r\xi\sigma]\right] + \Pr[\|a_i\| \geq \xi\sigma] \leq O((r\xi\sigma + \mu)\nu) + e^{-\Omega(\nu^2)} + e^{-\Omega(\xi^2 d)} \ .$$

If we choose $\xi = 1 + \Theta(\sqrt{d^{-1} \log \frac{1}{\alpha_0}})$ and $\nu = 2\kappa + \Theta(\sqrt{\log \frac{1}{\alpha_0}})$, then

$$\Pr\left[b_i \langle x, a_i \rangle \in [1 - 2r\xi\sigma - \mu, 1 + 2r\xi\sigma]\right] + \Pr[\|a_i\| \geq \xi\sigma] \leq O\left((r\sigma + \mu)\kappa \log \frac{1}{\alpha_0}\right) + \frac{\alpha_0}{5} \ .$$

$\square$

The next lemma extends Lemma 7.3 in two directions. First, it applies Chernoff bound to turn the probability upper bound in Lemma 7.3 into an upper bound on the actual number of data points; second, it applies standard epsilon-net argument to turn Lemma 7.3 into a "for all" argument with respect to all $x$ and to all $r$.

**Lemma 7.4.** *Under Assumption 3, for every $\alpha_0 \in (0,1)$ there exists some parameter $\xi \geq 1$ such that, as long as $n \geq \Omega(\frac{d}{\alpha_0} \log \frac{B\sigma\kappa}{\alpha_0})$, with probability at least $1 - e^{-\Omega(\alpha_0 n)}$, for every $x \in \mathbb{R}^d$ with $\|x\| \leq B$ and every $r \geq 0$,*

$$\left|\{i \in [n] \colon b_i \langle x, a_i \rangle \in [1 - r\xi\sigma - \mu, 1 + r\xi\sigma] \vee \|a_i\| \geq \xi\sigma\}\right| \leq n \cdot \left(O\left((r\sigma + \mu)\kappa \log \frac{1}{\alpha_0}\right) + \frac{\alpha_0}{2}\right) \ .$$

*Proof of Lemma 7.4.* We first apply Lemma 7.3 and Chernoff bound to derive that, for fixed $x \in \mathbb{R}^d$ and fixed $r \geq 0$, with probability at least $1 - e^{-\Omega(\alpha_0 n)}$ over the randomness of $a_1, \ldots, a_n$,

$$\left|\{i \in [n] \colon b_i \langle x, a_i \rangle \in [1 - 2r\xi\sigma - \mu, 1 + 2r\xi\sigma] \vee \|a_i\| \geq \xi\sigma\}\right| \leq n \cdot \left(O\left((r\sigma + \mu)\kappa \log \frac{1}{\alpha_0}\right) + \frac{\alpha_0}{4}\right) \ . \tag{7.1}$$

Next, we essentially want to take union bound with respect to all possible $x$ and $r$.

As for $r \geq 0$, it suffices for us to consider only $\Omega\left(\frac{\alpha_0}{\sigma\kappa \log(1/\alpha_0)}\right) \leq r \leq O\left(\frac{1}{\sigma\kappa \log(1/\alpha_0)}\right)$. This is because if $r$ is larger than this upper bound then the lemma becomes trivial, and if $r$ is smaller than this lower bound then we can replace it with this lower bound and then proceed the proof.

Furthermore, because we are hiding constants inside the big-$O$ notation, it suffices to consider only finitely many values $r$ in this interval: for instance, $r = r_0 \cdot 1.1^k$ for every $k \in \{0, 1, 2, \dots\}$ where $r_0 = \Theta\big(\frac{\alpha_0}{\sigma \kappa \log(1/\alpha_0)}\big)$.

As for $x$, we cover the space of $\|x\| \leq B$ by an $\varepsilon$-net with $\varepsilon = r_0/10$. This net has $e^{O(d \log(B/r_0))}$ many points $x$ and satisfies that each $x$ with $\|x\| \leq B$ is $\varepsilon$-close to at least one point in this net.

Applying union bound, we know that as long as $\alpha_0 n \geq \Omega(d \log \frac{B\sigma\kappa}{\alpha_0})$, the above (7.1) holds for every $x$ in this $\varepsilon$-net and every $r = r_0 1.1^k$ in the interval. It is not hard to derive that as a consequence, for all $x$ with $\|x\| \leq B$ and all $\Omega\big(\frac{\alpha_0}{\sigma\kappa\log(1/\alpha_0)}\big) \leq r \leq O\big(\frac{1}{\sigma\kappa\log(1/\alpha_0)}\big)$, it satisfies

$$\big|\{i \in [n] \colon b_i \langle x, a_i\rangle \in [1 - r\xi\sigma - \mu, 1 + r\xi\sigma] \vee \|a_i\| \geq \xi\sigma\}\big| \leq n \cdot \left(O\big((r\sigma + \mu)\kappa \log \frac{1}{\alpha_0}\big) + \frac{\alpha_0}{4}\right) . \tag{7.2}$$

(7.2) then implies for all $x$ with $\|x\| \leq B$ and all $r \geq 0$,

$$\big|\{i \in [n] \colon b_i \langle x, a_i\rangle \in [1 - r\xi\sigma - \mu, 1 + r\xi\sigma] \vee \|a_i\| \geq \xi\sigma\}\big| \leq n \cdot \left(O\big((r\sigma + \mu)\kappa \log \frac{1}{\alpha_0}\big) + \frac{\alpha_0}{2}\right) . \qquad \square$$

We are now ready to prove Theorem 7.2.

*Proof of Theorem 7.2.* We first observe that for each $i \in [n]$ and $x \in \mathbb{R}^d$,

$$b_i \langle x, a_i\rangle \notin [1 - r\xi\sigma - \mu, 1 + r\xi\sigma] \bigwedge \|a_i\| \leq \xi\sigma$$
$$\implies \forall y \text{ with } \|y - x\| \leq r \colon \nabla f_i(x) = \nabla f_i(y) .$$

This is because, if $b_i\langle x, a_i\rangle \leq 1 - r\xi\sigma - \mu$ and $\|a_i\| \leq \xi\sigma$ then $b_i\langle y, a_i\rangle \leq 1 - \mu$ so $\nabla f_i(x) = \nabla f_i(y) = -b_i a_i$, and similarly if $b_i\langle x, a_i\rangle \geq 1 + r\xi\sigma$ and $\|a_i\| \leq \xi\sigma$ then $\nabla f_i(x) = \nabla f_i(y) = b_i a_i$. As a result

$$\frac{|B(x,r)|}{n} = 1 - \frac{1}{n}\big|\{i \in [n] \colon \nabla f_i(x) = \nabla f_i(y) \text{ for all } y \text{ with } \|y - x\| \leq r\}\big|$$
$$\leq 1 - \frac{1}{n}\big|\{i \in [n] \colon b_i\langle x, a_i\rangle \notin [1 - r\xi\sigma - \mu, 1 + r\xi\sigma] \bigwedge \|a_i\| \leq \xi\sigma\}\big|$$
$$= \frac{1}{n}\big|\{i \in [n] \colon b_i\langle x, a_i\rangle \in [1 - r\xi\sigma - \mu, 1 + r\xi\sigma] \bigvee \|a_i\| > \xi\sigma\}\big|$$
$$\leq O\big((r\sigma + \mu)\kappa \log \frac{1}{\alpha_0}\big) + \frac{\alpha_0}{2} . \qquad \square$$

## Acknowledgements

We would like to thank Greg Yang, Ilya Razenshteyn and S ébastien Bubeck for discussing the motivation of this problem.

# Appendix

## A Appendix for Section 3

### A.1 Useful Lemma: Property of Index Sets

Recall from Definition 3.1 that if $(k_0, \ldots, k_t)$ is $k$'s lowbit sequence, then
$$\Lambda_k \stackrel{\text{def}}{=} \bigcup_{i=0}^{t-1} \left( B_{k_i}(k - k_i) \setminus B_{k_i}(k_{t-1} - k_i) \right) \quad \text{where} \quad B_k(r) \stackrel{\text{def}}{=} \Lambda_k \cap B(x_k, r \cdot \xi) \ .$$

For analysis purpose, we also define
$$\Lambda_{k,i} \stackrel{\text{def}}{=} B_{k_i}(k - k_i) \setminus B_{k_i}(k_{t-1} - k_i) \quad \text{and} \quad \Lambda_{k,i}^\perp \stackrel{\text{def}}{=} \Lambda_{k_i} \setminus B_{k_i}(k_t - k_i) \ .$$

Intuitively, we will show that $\Lambda_k$ contains the indices of $[n]$ where the stochastic gradients $\nabla f(x_k)$ need to be re-calculated. The remaining indices in $[n] \setminus \Lambda_k$ must have been stored in the memory.

Our main purpose of this subsection is to prove the following lemma:

**Lemma A.1.** *We have the following properties*

*(a)* $\Lambda_{k,i} = B_{k_i}(k_t - k_i) \setminus \left( \Lambda_{k_{i+1}} \cup \cdots \cup \Lambda_{k_{t-1}} \right)$

*(b)* $\Lambda_{k,i}^\perp = \left( \Lambda_{k_i} \setminus B_{k_i}(k_t - k_i) \right) \setminus \left( \Lambda_{k_{i+1}} \cup \cdots \cup \Lambda_{k_{t-1}} \right)$

*(c) The subsets* $\Lambda_{k,0}, \ldots, \Lambda_{k,t-1}$ *and* $\Lambda_{k,0}^\perp, \ldots, \Lambda_{k,t-1}^\perp$ *are pairwise disjoint.*

*(d)* $\bigcup_{i=0}^{t-1} \Lambda_{k,i} \cup \Lambda_{k,i}^\perp = [n]$.

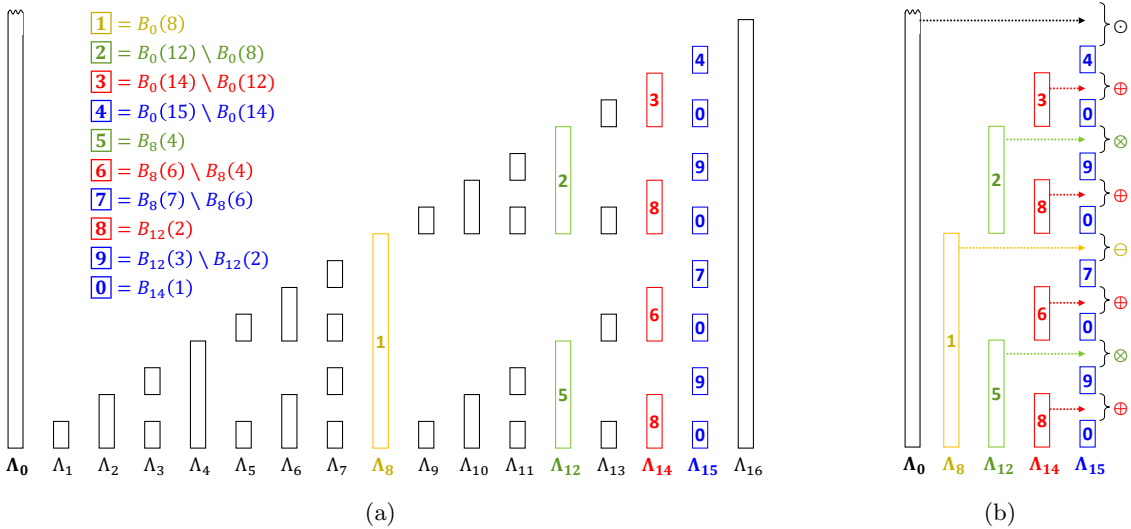

Figure 7: Illustration of index sets $\Lambda_k$ (same as Figure 2, repeated here for convenience)

Connecting back to Figure 2 (which we repeat here as Figure 7 for convenience):

- Lemma A.1.a says for instance box 7 equals $B_8(7) \setminus \Lambda_{12} \cup \Lambda_{14}$.
- Lemma A.1.b says for instance the $\ominus$ region does not intersect with $\Lambda_{12} \cup \Lambda_{14}$.
- Lemma A.1.c says for instance boxes 4, 0, 9, 7 and the $\ominus, \oplus, \odot, \otimes$ regions are pairwise disjoint.
- Lemma A.1.d says for instance boxes 4, 0, 9, 7 and the $\ominus, \oplus, \odot, \otimes$ regions jointly cover $[n]$.

We establish the following claim before proving Lemma A.1.

**Claim A.2.** *Suppose $(k_0, k_1, \ldots, k_b)$ is the lowbit sequence of $k = k_b$, then*
$$\forall a = 0, 1, \ldots, b-1: \quad \Lambda_{k_a} \bigcap \left( \Lambda_{k_{a+1}} \cup \cdots \cup \Lambda_{k_b} \right) = B_{k_a}(k_b - k_a)$$

For instance, in Figure 7, this means $\Lambda_8 \bigcap (\Lambda_{12} \cup \Lambda_{14}) = B_8(6)$.

*Proof of Claim A.2.* We prove by induction on $b$. If $b = 1$ then we know $(k_0, k_1)$ is the lowbit sequence of $k_1$, and thus by Definition 3.1 we have
$$\Lambda_{k_1} = B_{k_0}(k_1 - k_0) \setminus B_{k_0}(k_0 - k_0) = B_{k_0}(k_1 - k_0)$$
and thus
$$\Lambda_{k_0} \bigcap \Lambda_{k_1} = [n] \bigcap \Lambda_{k_1} = B_{k_0}(k_1 - k_0) \ .$$

Suppose the claim is true for $b - 1$, now we want to prove it for $b \geq 2$.

We first note that by induction (using the fact that $(k_0, k_1, \ldots, k_{b-1})$ is the lowbit sequence of $k_{b-1}$), we have
$$\forall i = 0, 1, \ldots, b-2: \quad \Lambda_{k_i} \bigcap \left( \Lambda_{k_{i+1}} \cup \cdots \cup \Lambda_{k_{b-1}} \right) = B_{k_i}(k_{b-1} - k_i) \ .$$

This means for every $i = 0, 1, \ldots, b-2$,
$$\Lambda_{k_b, i} = B_{k_i}(k_b - k_i) \setminus B_{k_i}(k_{b-1} - k_i) = B_{k_i}(k_b - k_i) \setminus \left( \Lambda_{k_{i+1}} \cup \cdots \cup \Lambda_{k_{b-1}} \right) \ . \quad \text{(A.1)}$$

We now prove the "$\supseteq$" direction of Claim A.2. Using (A.1), we have
$$\Lambda_{k_b} \supseteq \Lambda_{k_b, a} = B_{k_a}(k_b - k_a) \setminus B_{k_a}(k_{b-1} - k_a) = B_{k_a}(k_b - k_a) \setminus \left( \Lambda_{k_{a+1}} \cup \cdots \cup \Lambda_{k_{b-1}} \right) \ .$$

After rearranging, we have
$$\Lambda_{k_{a+1}} \cup \cdots \cup \Lambda_{k_b} \supseteq B_{k_a}(k_b - k_a) \ .$$

This proves the "$\supseteq$" direction because $B_{k_a}(k_b - k_a) \subseteq \Lambda_{k_a}$.

We next prove the "$\subseteq$" direction of Claim A.2. We have by induction
$$\Lambda_{k_a} \bigcap \left( \Lambda_{k_{a+1}} \cup \cdots \cup \Lambda_{k_{b-1}} \right) = B_{k_a}(k_{b-1} - k_a) \subseteq B_{k_a}(k_b - k_a)$$
and thus it suffices to show that
$$\Lambda_{k_a} \bigcap \left( \Lambda_{k_b} \setminus \left( \Lambda_{k_{a+1}} \cup \cdots \cup \Lambda_{k_{b-1}} \right) \right) \subseteq B_{k_a}(k_b - k_a) \ . \quad \text{(A.2)}$$

Recall $\Lambda_{k_b} = \Lambda_{k_b, 0} \cup \cdots \cup \Lambda_{k_b, b-1}$.

- For each $\Lambda_{k_b, i}$ with $i < a$, we have, owing to (A.1),
$$\Lambda_{k_b, i} = B_{k_i}(k_b - k_i) \setminus \left( \Lambda_{k_{i+1}} \cup \cdots \cup \Lambda_{k_{b-1}} \right)$$
so $\Lambda_{k_b, i}$ does not intersect with $\Lambda_{k_a}$.

- For each $\Lambda_{k_b, i}$ with $i > a$, we have
$$\Lambda_{k_b, i} = B_{k_i}(k_b - k_i) \setminus B_{k_i}(k_{b-1} - k_i) \subseteq \Lambda_{k_i},$$

To sum up, we have
$$\Lambda_{k_a} \bigcap \left( \Lambda_{k_b} \setminus \left( \Lambda_{k_{a+1}} \cup \cdots \cup \Lambda_{k_{b-1}} \right) \right) = \Lambda_{k_b, a} \subseteq B_{k_a}(k_b - k_a) \ .$$

This proves (A.2) and thus finishes the "$\subseteq$" direction of Claim A.2. □

*Proof of Lemma A.1.* These properties are corollaries of Claim A.2.

(a) Claim A.2 says we have $\Lambda_{k_i} \cap (\Lambda_{k_{i+1}} \cup \cdots \cup \Lambda_{k_{t-1}}) = B_{k_i}(k_{t-1} - k_i)$ and thus
$$\Lambda_{k,i} = B_{k_i}(k_t - k_i) \setminus B_{k_i}(k_{t-1} - k_i) = B_{k_i}(k_t - k_i) \setminus (\Lambda_{k_{i+1}} \cup \cdots \cup \Lambda_{k_{t-1}})$$

(b) Claim A.2 says we have $\Lambda_{k_i} \cap (\Lambda_{k_{i+1}} \cup \cdots \cup \Lambda_{k_{t-1}}) = B_{k_i}(k_{t-1} - k_i)$ and thus
$$\Lambda_{k,i}^\perp = \Lambda_{k_i} \setminus B_{k_i}(k_t - k_i) = (\Lambda_{k_i} \setminus B_{k_i}(k_t - k_i)) \setminus B_{k_i}(k_{t-1} - k_i)$$
$$= (\Lambda_{k_i} \setminus B_{k_i}(k_t - k_i)) \setminus (\Lambda_{k_{i+1}} \cup \cdots \cup \Lambda_{k_{t-1}})$$

(c) By Lemma A.1.a and Lemma A.1.b, we clearly have $\Lambda_{k,i}^\perp \cap \Lambda_{k,i} = \varnothing$.

For each $\Lambda_{k,i}$ and $i < j$, we have $\Lambda_{k,i} = B_{k_i}(k - k_i) \setminus (\Lambda_{k_{i+1}} \cup \cdots \cup \Lambda_{k_{b-1}})$ but also recall $\Lambda_{k,j} \subseteq \Lambda_{k_j}$ and $\Lambda_{k,j}^\perp \subseteq \Lambda_{k_j}$. Therefore we must have $\Lambda_{k,i} \cap \Lambda_{k,j} = \Lambda_{k,i} \cap \Lambda_{k,j}^\perp = \varnothing$.

For each $\Lambda_{k,i}^\perp$ and $i < j$, we have $\Lambda_{k,i}^\perp = (\Lambda_{k_i} \setminus B_{k_i}(k - k_i)) \setminus (\Lambda_{k_{i+1}} \cup \cdots \cup \Lambda_{k_{b-1}})$ but also recall $\Lambda_{k,j} \subseteq \Lambda_{k_j}$ and $\Lambda_{k,j}^\perp \subseteq \Lambda_{k_j}$. Therefore we must have $\Lambda_{k,i}^\perp \cap \Lambda_{k,j} = \Lambda_{k,i}^\perp \cap \Lambda_{k,j}^\perp = \varnothing$.

(d) By Lemma A.1.a and Lemma A.1.b, we have
$$\Lambda_{k,i} \cup \Lambda_{k,i}^\perp = \Lambda_{k_i} \setminus (\Lambda_{k_{i+1}} \cup \cdots \cup \Lambda_{k_{t-1}})$$
and therefore
$$\bigcup_{i=0}^{t-1} \Lambda_{k,i} \cup \Lambda_{k,i}^\perp = \bigcup_{i=0}^{t-1} \Lambda_{k_i} \setminus (\Lambda_{k_{i+1}} \cup \cdots \cup \Lambda_{k_{t-1}}) = \Lambda_{k_0} = [n] .$$ □

## A.2 Proof of Lemma 3.3

**Lemma 3.3.** *The construction of $\Lambda_0, \ldots, \Lambda_{m-1}$ ensures that $\mathbf{g} = \nabla f(x_k)$ in each iteration $k$.*

*Proof of Lemma 3.3.* Recall that in an iteration $k$, letting $(k_0, k_1, \ldots, k_t)$ be $k$'s lowbit sequence, then $[n] = \bigcup_{j=0}^{t-1} \Lambda_{k,j} \cup \Lambda_{k,j}^\perp$ according to Lemma A.1.d.

- For each $i \in \Lambda_k = \bigcup_{j=0}^{t-1} \Lambda_{k,j}$, we have recalculated $\mathbf{g}_i \leftarrow \nabla f_i(x_k)$ at Line 5 so it satisfies $\mathbf{g}_i = \nabla f_i(x_k)$.
- For each $i \in \Lambda_{k,j}^\perp = \Lambda_{k_j} \setminus B_{k_j}(k - k_j)$, we claim that $i$ never appeared in any $\Lambda_\ell$ for $\ell \in \{k_j + 1, k_j + 2, \ldots, k - 1\}$. To show this by induction on $\ell$.

  For each such $\ell$, first note that $k_j$ must also appear in the lowbit sequence of $\ell$. Let this lowbit sequence be $(\ell_0, \ell_1, \ldots, \ell_{t'})$ with $\ell = \ell_{t'}$. Suppose $\ell_{j'} = k_j$. Recall $\Lambda_\ell = \bigcup_{a=0}^{t'-1} \Lambda_{\ell,a}$.

  – For every $a = 0, 1, \ldots, j' - 1$, we have $\Lambda_{\ell,a} \cap \Lambda_{\ell_{j'}} = 0$ according to Lemma A.1.a but $i \in \Lambda_{k_j} = \Lambda_{\ell_{j'}}$, so $i \notin \Lambda_{\ell,a}$.
  – For every $a = j' + 1, \ldots, t' - 1$, we have $\Lambda_{\ell,a} \subseteq \Lambda_{\ell_a}$ but by induction $i \notin \Lambda_{\ell_a}$; therefore $i \notin \Lambda_{\ell,a}$.
  – For $a = j'$, we have $\Lambda_{\ell,a} = \Lambda_{\ell,j'} \subseteq B_{\ell_{j'}}(\ell - \ell_{j'}) \subseteq B_{\ell_{j'}}(k - \ell_{j'})$, but $k_j = \ell_{j'}$ and $i \in \Lambda_{k_j} \setminus B_{k_j}(k - k_j)$, so we must have $i \notin \Lambda_{\ell,a}$.

  In sum, we conclude that $i \notin \Lambda_\ell$ for any $\ell \in \{k_j + 1, k_j + 2, \ldots, k - 1\}$. In other words, the last time $\mathbf{g}_i$ was updated, it was at iteration $k_j$. This means $\mathbf{g}_i = \nabla f_i(x_{k_j})$. Using the facts that (1) $i \in \Lambda_{k_j}$, (2) $i \notin B_{k_j}(k - k_j) = \Lambda_{k_j} \cup B(x_{k_j}, (k - k_j) \cdot \xi)$, and (3) $\|x_k - x_{k_j}\| \leq (k - k_j) \cdot \xi$, we conclude that
  $$\mathbf{g}_i = \nabla f_i(x_{k_j}) = \nabla f_i(x_k) .$$

In sum, we have shown that at this iteration $k$, it satisfies $\mathbf{g}_i = \nabla f_i(x_k)$ for all $i \in [n]$. This means $\mathbf{g} = \nabla f(x_k)$. □

## A.3 Proof of Lemma 3.5

**Lemma 3.5.** *Under Assumption 2, we have $\frac{1}{n}\sum_{k=0}^{m-1}|\Lambda_k| \leq O(\alpha m + m^{1-\beta}\log^2 m)$ .*

*Proof of Lemma 3.5.* If $(k_0, \ldots, k_t)$ is the lowbit sequence of $k = k_t \geq 1$, then

$$|\Lambda_k| = \left|\bigcup_{i=0}^{t-1} B_{k_i}(k_t - k_i) \setminus B_{k_i}(k_{t-1} - k_i)\right| \leq \sum_{i=0}^{t-1}\left|B_{k_i}(k_t - k_i) \setminus B_{k_i}(k_{t-1} - k_i)\right| .$$

To upper bound $\sum_{k=1}^{m-1}|\Lambda_k|$, we group all terms $\left|B_{k_i}(k_t - k_i) \setminus B_{k_i}(k_{t-1} - k_i)\right|$ together with respect the same $j = k_i$.

For each $j = 0, 1, 2, \ldots, m-1$, suppose $\mathsf{lowbit}(j) = 2^p$ for $p \geq 0$. Then, it is easy to verify that $j$ appears in the lowbit sequence $(k_0, k_1, \ldots, k_t)$ of $k = k_t$ (with $j = k_i$ for some $i < t$) if and only if $k \in [j+1, j+2^p-1]$. For each $k$ in this range, we have

$$B_j(k_t - j) \setminus B_j(k_{t-1} - j) = B_j(k - j) \setminus B_j(k - \mathsf{lowbit}(k) - j)$$

and therefore we want to upper bound the cardinality of the right hand side over all $k \in [j+1, j+2^p-1]$. We calculate that:

$$\sum_{k=j+1}^{j+2^p-1}\left|B_j(k-j) \setminus B_j(k - \mathsf{lowbit}(k) - j)\right|$$

$$\stackrel{①}{=} \sum_{k=j+1}^{j+2^p-1}\left|B_j(k-j) \setminus B_j(k - \mathsf{lowbit}(k-j) - j)\right|$$

$$= \sum_{k=1}^{2^p-1}\left|B_j(k) \setminus B_j(k - \mathsf{lowbit}(k))\right|$$

$$\stackrel{②}{=} \sum_{a=0}^{p-1}\sum_{b=0}^{2^{p-a-1}-1}\left|B_j((2b+1)\cdot 2^a) \setminus B_j(2b \cdot 2^a)\right|$$

$$\stackrel{③}{=} \sum_{a=0}^{p-1}\left|\bigcup_{b=0}^{2^{p-a-1}-1}\left(B_j((2b+1)\cdot 2^a) \setminus B_j(2b \cdot 2^a)\right)\right|$$

$$\leq \sum_{a=0}^{p-1}\left|B_j(2^p)\right| \stackrel{④}{\leq} p \cdot \psi(2^p \xi) .$$

Above, equality ① is because $\mathsf{lowbit}(k-j) = \mathsf{lowbit}(k)$; equality ② is by re-grouping the indices $k \in \{1, 2, \ldots, 2^p - 1\}$ according to their lowest binary bit; equality ③ is because the subsets over $b = 0, 1, \ldots, 2^{p-a-1} - 1$ are pairwise disjoint; and inequality ④ is by our Assumption 2.

Finally, if $p_j = \mathsf{lowbit}(j)$ for each $j = 0, 1, \ldots, m-1$, then we have

$$\frac{1}{n}\sum_{k=1}^{m-1}|\Lambda_k| \leq \sum_{j=0}^{m-1} p_j \cdot \psi(2^{p_j}\xi) \leq \sum_{p=0}^{\lfloor\log_2 m\rfloor} p \cdot \psi(2^p \xi) \times \frac{m}{2^p} .$$

With our choice of $\psi(\cdot)$ and $\xi = C/m$, it satisfies

$$\sum_{p=0}^{\lfloor\log_2 m\rfloor} p \cdot \psi(2^p \xi) \times \frac{m}{2^p} \leq \sum_{p=0}^{\lfloor\log_2 m\rfloor} p \cdot (\alpha + (2^p \xi/C)^\beta) \times \frac{m}{2^p} = O(\alpha m + m^{1-\beta}\log^2 m) . \qquad \square$$

## A.4 Useful Lemma: Truncated Gradient Descent

We show a simple lemma on the performance of truncated gradient descent. Its proof is a simple modification of Nesterov [27] for proving the convergence of gradient descent.

**Lemma A.3.** *Given $x_0 \in \mathbb{R}^d$ and $\xi \geq 0$, suppose $f(x)$ is $L$-smooth, convex and we apply*

$$x_{k+1} \leftarrow x_k - \min\left\{\frac{\xi}{\|\nabla f(x_k)\|}, \frac{1}{L}\right\} \cdot \nabla f(x_k)$$

*for $m$ steps, then it satisfies*

$$\|x_m - x^*\| \leq \|x_0 - x^*\| \quad \text{and} \quad f(x_m) - f(x^*) \leq \max\left\{\frac{4L \cdot \|x_0 - x^*\|^2}{m}, f(x_0) - f(x^*) - \frac{mL\xi^2}{4}\right\} .$$

The classical $\frac{1}{T}$ convergence of gradient descent is a special case of Lemma A.3 by setting $\xi = +\infty$.

*Proof of Lemma A.3.* We first prove $\|x_m - x^*\| \leq \|x_0 - x^*\|$. We have

$$\|x_{k+1} - x^*\|^2 = \|x_k - x^*\|^2 - 2\langle x_k - x_{k+1}, x_k - x^*\rangle + \|x_k - x_{k+1}\|^2$$

and using smoothness (c.f. [27, 2.1.8]), we have $\langle \nabla f(x_k), x_k - x^*\rangle \geq \frac{1}{L}\|\nabla f(x_{k+1})\|^2$. This means,

1. if $x_{k+1} = x_k - \frac{1}{L}\nabla f(x_k)$, then it satisfies $\|x_{k+1} - x^*\|^2 \leq \|x_k - x^*\|^2 - \frac{1}{L^2}\|\nabla f(x_{k+1})\|^2$.

2. if $x_{k+1} = x_k - \frac{\xi}{\|\nabla f(x_k)\|}\nabla f(x_k)$, then it satisfies $\|x_{k+1} - x^*\|^2 \leq \|x_k - x^*\|^2 - \frac{2\xi}{L}\|\nabla f(x_{k+1})\| + \xi^2 \leq \|x_k - x^*\|^2 - \xi^2$.

In both cases, we have $\|x_{k+1} - x^*\| \leq \|x_k - x^*\| \leq \cdots \|x_0 - x^*\|$.

We next prove the upper bound on $f(x_m) - f(x^*)$. By the $L$-smoothness of $f(x)$, we have

$$f(x_k) - f(x_{k+1}) \geq \langle \nabla f(x_k), x_k - x_{k+1}\rangle - \frac{L}{2}\|x_k - x_{k+1}\|^2 .$$

Therefore,

1. if $x_{k+1} = x_k - \frac{1}{L}\nabla f(x_k)$, then it satisfies $f(x_k) - f(x_{k+1}) \geq \frac{1}{2L}\|\nabla f(x_k)\|^2$.

2. if $x_{k+1} = x_k - \frac{\xi}{\|\nabla f(x_k)\|}\nabla f(x_k)$, then it satisfies $f(x_k) - f(x_{k+1}) \geq \frac{\xi}{2}\|\nabla f(x_k)\| \geq \frac{L\xi^2}{2}$.

There are two possibilities. If case 2 above takes place for at least $m/2$ steps, then we have $f(x_0) - f(x_m) \geq \frac{m}{2} \cdot \frac{L\xi^2}{2}$ and therefore $f(x_m) - f(x^*) \leq f(x_0) - f(x^*) - \frac{mL\xi^2}{4}$. This finishes the proof. For the remainder of the proof, we assume that case 1 above takes place for at least $m/2$ steps.

Denote by $\Delta_k = f(x_k) - f(x^*)$. By convexity, we have

$$\Delta_k \leq \langle \nabla f(x_k), x_k - x^*\rangle \leq \|\nabla f(x_k)\| \cdot \|x_k - x^*\| \leq \|\nabla f(x_k)\| \cdot \|x_0 - x^*\| .$$

At the same time, we have $\Delta_{k+1} \leq \Delta_k$ and whenever case 1 happens, we have

$$\Delta_{k+1} \leq \Delta_k - \frac{1}{2L}\|\nabla f(x_k)\|^2 \leq \Delta_k - \frac{1}{2L}\frac{\Delta_k^2}{\|x_0 - x^*\|^2} .$$

Dividing both sides by $\Delta_k \cdot \Delta_{k+1}$, and using $\Delta_k \geq \Delta_{k+1}$, we conclude that

$$\frac{1}{\Delta_k} \leq \frac{1}{\Delta_{k+1}} - \frac{1}{2L \cdot \|x_0 - x^*\|^2} .$$

In other words, after $m$ iterations, since case 1 happens at least $m/2$ times, we have

$$\frac{1}{\Delta_m} \geq \frac{m}{4L \cdot \|x_0 - x^*\|^2} . \qquad \square$$

## A.5 Proof of Theorem 3.6

**Theorem 3.6.** *Given any $x^{(0)} \in \mathbb{R}^d$ and $D > 0$ that is an upper bound on $\|x^{(0)} - x^*\|^2$. Suppose Assumption 1 and 2 are satisfied with parameters $C \in (0, D], \alpha \in [0, 1], \beta \in (0, 1]$. Then, denoting by $m_s = \lceil (1 + \frac{C^2}{16D^2})^s \rceil$, we have that $\mathtt{GD^{lin}}(f, x_0, S, C, D)$ outputs a point $x \in \mathbb{R}^d$ satisfying $f(x) - f(x^*) \leq \frac{4LD^2}{m_S}$ with gradient complexity $T_{\mathsf{time}} = O\big(\sum_{s=1}^{S} \alpha m_s + m_s^{1-\beta} \log^2 m_s \big)$.*

*Proof of Theorem 3.6.* We first inductively prove that

$$\|x^{(s)} - x^*\| \leq D \quad \text{and} \quad f(x^{(s)}) - f(x^*) \leq \frac{4LD^2}{m_s} \ . \tag{A.3}$$

In the base case $s = 0$, recall $x^{(0)}$ comes from the initial vector $x^{(0)}$ after applying $m_0$ steps of gradient descent. By Lemma A.3 with $\xi = +\infty$, we derive that

$$\|x^{(0)} - x^*\| \leq \|x^{(0)} - x^*\| \leq D \quad \text{and} \quad f(x^{(0)}) - f(x^*) \leq \frac{4L\|x^{(0)} - x^*\|^2}{m_0} \leq \frac{4LD^2}{m_0} \ .$$

Suppose (A.3) holds for $s - 1$. In epoch $s \geq 1$, $\mathtt{GD^{lin}}$ applies truncated gradient with $\xi_s = \frac{C}{m_s}$, starting from $x_0 = x^{(s-1)}$, for $m_s$ steps. Using Lemma A.3, we have

$$\|x_{m_s} - x^*\| \leq \|x_0 - x^*\| \leq D \quad \text{and} \quad f(x_{m_s}) - f(x^*) \leq \max\left\{\frac{4L \cdot \|x_0 - x^*\|^2}{m_s}, \frac{4LD^2}{m_{s-1}} - \frac{m_s L \xi_s^2}{4}\right\} \ .$$

Since $\frac{4LD^2}{m_{s-1}} - \frac{m_s L \xi_s^2}{4} = \frac{4LD^2}{m_{s-1}} - \frac{LC^2}{4m_s} \geq \frac{4LD^2}{m_s}$ according to our choice of $m_s \leq (1 + \frac{C^2}{16D^2})m_{s-1}$, we conclude $f(x_{m_s}) - f(x^*) \leq \frac{4LD^2}{m_s}$ as desired. This proves (A.3).

As for the gradient complexity, Claim 3.4 tells us that it equals $\frac{1}{n} \sum_{k=0}^{m-1} |\Lambda_k|$. By Lemma 3.5, we have

$$T_{\mathsf{time}} = O\Big(\sum_{s=1}^{S} \big(\alpha m_s + m_s^{1-\beta} \log^2 m_s\big)\Big) \ . \qquad \square$$

## A.6 Proof of Theorem 3.7

**Theorem 3.7.** *In the setting of Theorem 3.6, given any $T \geq 1$, one can choose $S$ so that*

- *If $\beta = 1$, then $\mathtt{GD^{lin}}$ finds a point $x$ in gradient complexity $T_{\mathsf{time}} = O(T)$ s.t.*

$$f(x) - f(x^*) \leq O\big(\tfrac{LD^4}{C^2} \cdot \tfrac{\alpha}{T}\big) + \frac{LD^2}{2^{\Omega(C^2T/D^2)^{1/3}}} \ .$$

- *If $\beta \in (0, 1)$ is constant, $\mathtt{GD^{lin}}$ finds a point $x$ in gradient complexity $T_{\mathsf{time}} = O(T \log^2 T)$ s.t.*

$$f(x) - f(x^*) \leq O\big(\tfrac{LD^4}{C^2} \cdot \tfrac{\alpha}{T} + \tfrac{LD^{2+\frac{2}{1-\beta}}}{C^{\frac{2}{1-\beta}}} \cdot \tfrac{1}{T^{1/(1-\beta)}}\big) \ .$$

*Proof of Theorem 3.7.* Without loss of generality we assume $T \geq \Omega(D^2/C^2)$ because otherwise $f(x) - f(x^*) \leq f(x^{(0)}) - f(x^*) \leq O(LD^2)$ so the statements are trivially satisfied.

If $\beta = 1$, we choose $S$ large enough so that $m_S = \min\big\{\Theta(\frac{C^2 T}{D^2 \alpha}), 2^{\Theta(C^2T/D^2)^{1/3}}\big\}$. In this case,

$$T_{\mathsf{time}} = O\Big(\sum_{s=1}^{S} \alpha m_s + \log^2 m_s\Big) \leq O\Big(\alpha \frac{D^2}{C^2} m_S + S \cdot \log^2 m_S\Big) \leq O\Big(T + \frac{D^2}{C^2} \log^3 m_S\Big) \leq O(T)$$

and

$$f(x) - f(x^*) \leq O\Big(\frac{LD^2}{m_S}\Big) \leq O\Big(\frac{LD^4}{C^2} \cdot \frac{\alpha}{T}\Big) + \frac{LD^2}{2^{\Omega(C^2T/D^2)^{1/3}}} \ .$$

If $\beta \in (0,1)$ is a constant, we choose $S$ large enough so that $m_S = \Theta(\min\{\frac{C^2T}{D^2\alpha}), (\frac{TC^2}{D^2})^{1/(1-\beta)}\}$. In this case, we have

$$T_{\text{time}} = O\Big(\sum_{s=1}^{S} \alpha m_s + m_s^{1-\beta}\log^2 m_s\Big) \leq O\Big(\frac{D^2}{C^2}\alpha m_s + \frac{D^2}{C^2}m_S^{1-\beta}\log^2 m_S\Big) \leq O(T\log^2 T)$$

and

$$f(x) - f(x^*) \leq O\Big(\frac{LD^2}{m_S}\Big) \leq O\Big(\frac{LD^4}{C^2}\cdot\frac{\alpha}{T} + \frac{LD^{2+\frac{2}{1-\beta}}}{C^{\frac{2}{1-\beta}}}\cdot\frac{1}{T^{1/(1-\beta)}}\Big) \ . \qquad \Box$$

## B  Implementation Details

In this section we discuss some implementation details regarding `SVRG`$^{\text{lin}}$, `SCSG`$^{\text{lin}}$ and `GD`$^{\text{lin}}$.

### B.1  Implementation Details for `SVRG`$^{\text{lin}}$

We discuss the implementation details of `SVRG`$^{\text{lin}}$ (pseudocode in Algorithm 2) in this section. To efficiently implement `SVRG`$^{\text{lin}}$, we need the following special treatments:

- The algorithm always maintains $\{H_s\}_s$ as well as $[n] \setminus (H_0 \cup \cdots \cup H_s)$.
- The algorithm always maintains $\sum_{i \in H_s} \nabla f_i(x^{(s)})$ whenever an index moves out of $H_s$.
- At the start of epoch $s$, the full gradient $\nabla f(x_0)$ should be computed in gradient complexity $|H_s|/n$ in Line 3 of `SVRG`$^{\text{lin}}$. This is because, for any past epoch $s' < s$, the value $\sum_{i \in H_{s'}} \nabla f_i(x^{(s')})$ can be retrieved in $O(1)$ time from storage. We only need to compute gradients for indices $i \in H_s = [n] \setminus (H_1 \cup \cdots \cup H_{s-1})$.
- For each iteration in an epoch, since we have maintained the set $[n] \setminus (H_0 \cup \cdots \cup H_s)$, we can sample $i$ in $O(1)$ time. As for gradient complexity, `SVRG`$^{\text{lin}}$ only needs to compute new gradient $\nabla f_i(x_k)$ once, with gradient complexity $1/n$.
- After $H_s$ is constructed, for future use, we should calculate the lingering radius $\delta(x^{(s)}, i)$ for all $i \in H_s$. According to Assumption 1, this requires the same gradient complexity $|H_s|/n$ as computing $\nabla f(x_0)$. Then, we sort all the indices $i \in H_s$ in their increasing order of $\delta(x^{(s)}, i)$.
- In Line 11 of `SVRG`$^{\text{lin}}$, for each $s' \leq s$, thanks to the fact that indices in $H_{s'}$ are sorted, we can efficiently determine all the indices $i$ from $H_{s'}$ with $\delta_i(x^{(s')}, i) < r$ for any given parameter $r$. Therefore, we can remove each of such index in $O(1)$ time. This does not require any new computation of $\delta(x^{(s')}, i)$, and thus does not affect gradient complexity.

The above treatments ensure that the total gradient complexity of an epoch of length $m$ is $O(|H_s|/n) = O(m/n)$. Therefore, the amortized per-iteration gradient complexity of `SVRG`$^{\text{lin}}$ is $O(1)$, the same as SVRG. On the other hand, because `SVRG`$^{\text{lin}}$ has reduced variance of the gradient estimator **g** comparing to SVRG (as we discussed in Section 4), it should converge in a smaller number of iterations.

*Remark* B.1. Above, we did not discuss the additional time complexity (in addition to gradient complexity) for `SVRG`$^{\text{lin}}$. Most of the operations require $O(d)$ additional time per iteration, except two of them:

- To perform sorting on $H_s$ at the beginning of an epoch, we need $O(m\log m)$ total time for each epoch, amortized to $O(\log m)$ per iteration so is usually negligible.

- To calculate $\|x^{(s')} - x_{k+1}\|$ for every $s' \leq s$ in Line 11 of SVRG$^{\mathrm{lin}}$, we need $O(Sd)$ time per iteration. This can be bottleneck for certain applications.[12] If this happens, we propose a simple trick to reduce this complexity to $O(d)$:

    - Within an epoch, we calculate $\{\|x^{(s')} - x_{k+1}\|\}_{s'=0}^{s}$ exactly only once every $S$ iterations.
    - For the remaining iterations, we use triangle inequality to give upper bound to $\|x^{(s')} - x_{k+1}\|$. For instance, if we have calculated $\|x^{(s')} - x_{k'}\|$ exactly, then we can use $\|x^{(s')} - x_{k'}\| + \|x_{k'} - x_{k+1}\|$ to upper bound $\|x^{(s')} - x_{k'}\|$.

    We demonstrate in Section 5 that this trick works well in practice.

## B.2 Implementation Details for SCSG$^{\mathrm{lin}}$

For SCSG$^{\mathrm{lin}}$, we set the length of the initial epoch to be $\overline{m}_0 = 100$. Other practical treatments for SCSG$^{\mathrm{lin}}$ are the same as those for SVRG$^{\mathrm{lin}}$.

## B.3 Implementation Details for GD$^{\mathrm{lin}}$

To improve the practical performance of GD$^{\mathrm{lin}}$, we make several natural changes to Algorithm 1:

- Before starting GD$^{\mathrm{lin}}$, we run GD for 50 steps, so that the initial solution for GD$^{\mathrm{lin}}$ has a more significant lingering effect. The running time of GD$^{\mathrm{lin}}$ reported in Section 6 has included these 50 steps.
- We set $m = \min\{1000, \lceil 100 \times 1.1^s \rceil\}$, where the base 1.1 is selected for empirical performance, and the upper bound 1000 is selected to control memory usage.
- We view $C$ as the *learning rate*, and replace Line 8 of Algorithm 1 with $x_{k+1} \leftarrow x_k - \frac{\xi}{\|\mathbf{g}\|}\mathbf{g}$. In other words, the distance traveled by GD$^{\mathrm{lin}}$ in each epoch becomes exactly $C$.

Below, we discuss how to efficiently compute the index sets $\Lambda_k$:

- Each time a set $\Lambda_k$ is constructed, we sort all the indices in $\Lambda_k$ in increasing order of their lingering radii. The running time of this step is $O(|\Lambda_k|d)$ for computing the lingering radii, plus $O(|\Lambda_k| \log |\Lambda_k|)$ for sorting.
- Then, every $B_k(r)$ can be represented by a single integer $l$ such that $B_k(r)$ is the set of $l$ indices with the smallest lingering radii in $\Lambda_k$. We use a hash table to store all the values of $l$ as a function of $(k, r)$.
- To compute $B_{k_i}(k - k_i) \setminus B_{k_i}(k_{t-1} - k_i)$, we perform the following steps:
    1. Look in the hash table for $l = |B_{k_i}(k_{t-1} - k_i)|$.
    2. Starting from $l$, we find an $l' \geq l$ such that $l' = |B_{k_i}(k - k_i)|$.
    3. Store $l' = |B_{k_i}(k - k_i)|$ in the hast table.

    The time complexity of these steps is $O(|B_{k_i}(k - k_i)| - |B_{k_i}(k_{t-1} - k_i)|)$.

---
[12]Indeed, if computing a stochastic gradient $\nabla f_i(x)$ can be done in time $O(d)$, calculating those norms becomes the bottleneck for each iteration.